\documentclass[a4paper]{easychair}

\usepackage{bussproofs}
\usepackage{amsmath,amssymb}
\usepackage{mathrsfs}
\usepackage{enumerate}
\usepackage{stmaryrd}
\usepackage{tikz}
\tikzset{shorten <>/.style={shorten >=#1, shorten <=#1}}
\usetikzlibrary{patterns,arrows}
\usetikzlibrary{automata, positioning}
\usepackage{setspace}
\newtheorem{theorem}{Theorem}[section]
\newtheorem{lemma}[theorem]{Lemma}
\newtheorem{proposition}[theorem]{Proposition}
\newtheorem{corollary}[theorem]{Corollary}
\newtheorem{fact}[theorem]{Fact}
\newtheorem{remark}[theorem]{Remark}
\theoremstyle{plain}
\newtheorem{definition}[theorem]{Definition}
\newtheorem{example}[theorem]{Example}
\renewcommand{\phi}{\varphi}

\newcommand{\bd}{\blacklozenge}
\newcommand{\bb}{\blacksquare}
\newcommand{\sub}{\subseteq}

\newcommand{\comb}[1]{(#1)^{c}}

\newcommand{\ve}{\varnothing}
\newcommand{\D}{\Diamond}
\newcommand{\B}{\Box}

\newcommand{\cset}[1]{{\{ #1 \}}}
\newcommand{\tup}[1]{{\langle #1 \rangle}}

\newcommand{\rsto}{{\upharpoonright}}
\newcommand{\LT}{\mathsf{NExt}(\mathsf{K4}_t)}
\newcommand{\K}{\mathsf{NExt}(\mathsf{K}_t)}
\newcommand{\NExt}{\mathsf{NExt}}
\newcommand{\iso}{\cong}

\newcommand{\Fr}{\mathsf{Fr}}

\newcommand{\R}{R_\sharp}

\newcommand{\gf}{\mathbb{F}}
\renewcommand{\gg}{\mathbb{G}}
\newcommand{\F}{\mathfrak{F}}

\newcommand{\G}{\mathfrak{G}}

\newcommand{\M}{\mathfrak{M}}
\newcommand{\C}{\mathfrak{Ch}}
\renewcommand{\L}{ \mathscr{L}_t}
\newcommand{\ST}{{\mathsf{S4}_t}}

\newcommand{\J}{\mathcal{J}}

\newcommand{\md}{\models}

\usepackage{tikz}
\usetikzlibrary{patterns,arrows}

\title{
    Degree of Kripke-incompleteness of Tense Logics
    }

\author{
    Qian Chen
}

\institute{
    The Tsinghua-UvA JRC for Logic, Department of Philosophy, Tsinghua University, China  \and  
    Institute for Logic, Language and Computation, University of Amsterdam,  The Netherlands
}

\authorrunning{~}

\titlerunning{~}

\begin{document}

\maketitle

\begin{abstract}
    The degree of Kripke-incompleteness of a logic $L$ in some lattice $\mathcal{L}$ of logics is the cardinality of logics in $\mathcal{L}$ which share the same class of Kripke-frames with $L$. A celebrated result on Kripke-incompleteness is Blok's dichotomy theorem for the degree of Kripke-incompleteness in $\mathsf{NExt}(\mathsf{K})$: every modal logic $L\in\mathsf{NExt}(\mathsf{K})$ is of the degree of Kripke-incompleteness $1$ or $2^{\aleph_0}$. In this work, we show that the dichotomy theorem for $\mathsf{NExt}(\mathsf{K})$ can be generalized to the lattices $\K$, $\LT$ and $\NExt(\ST)$ of tense logics. We also prove that in $\K$, $\LT$ and $\NExt(\ST)$, iterated splittings are exactly the strictly Kripke-complete logics.%
\end{abstract}

\section{Introduction}

A logic $L$ is Kripke-complete if $L$ is the logic of some class of Kripke-frames. Kripke-completeness of modal logics has been extensively studied since 1960s. Thomason \cite{Thomason1972} established the existence of Kripke-incomplete tense logics. Fine~\cite{Fine1974a} and van Benthem~\cite{VanBenthem1978} gave examples of Kripke-incomplete modal logics. To study Kripke-completeness at a higher level, Fine \cite{Fine1974a} introduced the degree of Kripke-incompleteness of logics. Let $\Fr(L)$ denote the class of all frames validating $L$. For any lattice $\mathcal{L}$ of logics and $L\in\mathcal{L}$, the {\em degree of Kripke-incompleteness $\mathsf{deg}_{\mathcal{L}}(L)$ of $L$ in $\mathcal{L}$} is defined as:
\begin{center}
    $\mathsf{deg}_{\mathcal{L}}(L)=|\cset{L'\in\mathcal{L}:\Fr(L')=\Fr(L)}|$.
\end{center}
In general, studying the degree of Kripke-incompleteness in $\mathcal{L}$ amounts to analyzing the equivalence relation $\equiv_\Fr$ on $\mathcal{L}$, where $L_1\equiv_\Fr L_2$ iff $L_1$ shares the same class of frames as $L_2$, i.e., $\Fr(L_1)=\Fr(L_2)$. The degree of Kripke-incompleteness of $L$ is the cardinality of the equivalence class $[L]_{\equiv_\Fr}$ in $\mathcal{L}$. A logic $L$ is \emph{strictly Kripke-complete in $\mathcal{L}$} if $\mathsf{deg}_\mathcal{L}(L)=1$.

A celebrated result in this field is Blok's dichotomy theorem for the degree of Kripke-incompleteness in $\mathsf{NExt}(\mathsf{K})$: every modal logic $L\in\mathsf{NExt}(\mathsf{K})$ is of the degree of Kripke-incompleteness $1$ or $2^{\aleph_0}$. This theorem was first proved in \cite{Blok1978a} algebraically by showing that union-splittings in $\mathsf{NExt}(\mathsf{K})$ are exactly the consistent strictly Kripke-complete logics and all other consistent logics have the degree $2^{\aleph_0}$. Blok's characterization shows the connection between strictly Kripke-complete normal modal logics and splittings of lattices of logics in $\NExt(\mathsf{K})$. For more research on splittings of lattices of modal, tense and subframe logics, we refer the readers to \cite{Wolter1997c,Kracht1993,Kracht1990,Rautenberg1980}.
A proof based on relational semantics was given later in \cite[Section~10.5]{Chagrov.Zakharyaschev1997}. The characterization of the degree of Kripke-incompleteness given by Blok indicates locations of Kripke-complete logics in the lattice $\mathsf{NExt}(\mathsf{K})$. Since Blok's proof relies heavily on non-transitive frames, it is natural to ask whether the dichotomy theorem holds for sublattices of $\NExt(\mathsf{K})$, especially for the lattices of transitive modal logics such as $\mathsf{K4}$ and $\mathsf{S4}$. These problems remain open, see \cite[Problem~10.5]{Chagrov.Zakharyaschev1997}.

Generally, one can always replace the class $\Fr$ of all Kripke frames with some proper class $\mathcal{C}$ of mathematical structures, for example, the class $\mathsf{MA}$ of all modal algebras, the class $\mathsf{NF}$ of all neighborhood frames and the class $\mathsf{Fin}$ of all finite frames. Let $\mathcal{L}=\NExt(\mathsf{K})$. Since every normal modal logic is complete with respect to modal algebras \cite[Theorem~7.73]{Chagrov.Zakharyaschev1997} we obtain that $\equiv_\mathsf{MA}$ is the identity relation on $\mathcal{L}$. Note that logics that enjoy the finite model property (FMP) are Kripke complete and Kripke complete logics are neighborhood complete, we have $\equiv_\mathsf{NF}\sub\equiv_\mathsf{Fr}\sub\equiv_\mathsf{Fin}$. 
The degree of modal incompleteness with respect to neighborhood semantics was also well-investigated, e.g., by Chagrova \cite{Chagrova1998}, Dziobiak \cite{Dziobiak1978} and Litak \cite{Litak2004}. Dziobiak \cite{Dziobiak1978} proved the dichotomy theorem for the degree of incompleteness in the lattice $\mathsf{Ext}(\mathsf{D}\oplus(\B^np\to\B^{n+1}p))$ w.r.t neighborhood semantics for all $n\in\omega$. Litak \cite{Litak2004} studied modal incompleteness w.r.t Boolean algebras with operators (BAOs) and showed the existence of a continuum of neighborhood-incomplete modal logics extending $\mathsf{Grz}$. For more on modal incompleteness from an algebraic view, we refer the readers to \cite{Litak2005}.
Bezhanishvili et al.~\cite{Bezhanishvili.Bezhanishvili.ea2025} introduced the notion of the degree of FMP of $L$ in $\mathcal{L}$, which is in fact the cardinality of the equivalence class $[L]_{\equiv_\mathsf{Fin}}$. The anti-dichotomy theorem for the degree of FMP for extensions of the intuitionistic propositional logic $\mathsf{IPC}$ was proved in \cite{Bezhanishvili.Bezhanishvili.ea2025}: for each cardinal $\kappa$ with $0<\kappa\leq\aleph_0$ or $\kappa=2^{\aleph_0}$, there exists $L\in\mathsf{Ext}(\mathsf{IPC})$ such that the degree of FMP of $L$ in $\mathsf{Ext}(\mathsf{IPC})$ is $\kappa$. It was also shown in \cite{Bezhanishvili.Bezhanishvili.ea2025} that the anti-dichotomy theorem of the degree of FMP holds for $\mathsf{NExt}(\mathsf{K4})$ and $\mathsf{NExt}(\mathsf{S4})$. Degrees of FMP in lattices of bi-intuitionistic logics were studied in \cite{Chernev2022}. Given close connections between modal logics, bi-intuitionistic logics and tense logics, it is natural to study the degree of Kripke-incompleteness in lattices of tense logics.

Tense logics are bi-modal logics that include a future-looking necessity modality $\B$ and a past-looking possibility modality $\bd$, of which the lattices are substantially different from those of modal logics (see \cite{Kracht1992,Thomason1972,Ma.Chen2021}). As far as we are aware, the degree of Kripke-incompleteness in lattices of tense logics has not been investigated systematically. In this work, we study Kripke-incompleteness in lattices of tense logics. We start with the lattice $\K$ of all tense logics. Inspired by the proof for Blok's dichotomy theorem in \cite{Chagrov.Zakharyaschev1997}, we prove the dichotomy theorem for tense logics, that is, every tense logic $L\in\K$ is of degree of Kripke-incompleteness $1$ or $2^{\aleph_0}$. This is proved by showing that union-splittings in $\K$ are exactly the strictly Kripke-complete logics and all other logics have the degree $2^{\aleph_0}$. By a similar argument, we prove the dichotomy theorem of the degree of Kripke-incompleteness for $\LT$. Finally, we turn to the lattice $\NExt(\ST)$. We provide the following characterization of the degree of Kripke-incompleteness in $\NExt(\ST)$: iterated splittings are strictly Kripke-complete and all other logics are of degree $2^{\aleph_0}$, where iterated splittings are intuitively splittings in the lattice of splitting logics. The dichotomy theorem of the degree of Kripke-incompleteness for $\NExt(\ST)$ follows from the characterization immediately. It also follows that in the lattice $\NExt(\ST)$, strictly Kripke-complete logics are no longer the union-splittings.

Results obtained in this work indicate that the notion of iterated splitting fits better with strict Kripke-completeness. Wolter \cite{Wolter1997c} studied the iterated splittings of tense and subframe logics. Blok \cite{Blok1978a} showed that consistent iterated splittings in $\NExt(\mathsf{K})$ are exactly the union-splittings. In this work, we show that in the lattices $\NExt(\mathsf{K}_t)$ and $\NExt(\mathsf{K4}_t)$, iterated splittings are exactly union-splittings. Hence, we obtain the following unified characterization of the degree of Kripke-incompleteness in the lattices $\K$, $\LT$ and $\NExt(\ST)$: iterated splittings are strictly Kripke-complete and all other logics are of degree $2^{\aleph_0}$.

This paper is structured as follows: Section~2 gives preliminaries on tense logics, splittings and the degree of Kripke-incompleteness. Section~3 introduces reflective unfolding of Kripke frames, which is one of the most important method used in this paper. Sections~4 and 5 prove the dichotomy theorem of the degree of Kripke-incompleteness for $\K$, $\LT$ and $\NExt(\ST)$. Section~6 gives some concluding remarks.

\section{Preliminaries}

Basic notations on modal and tense logic can be found in e.g. \cite{Blackburn.deRijke.ea2001,Chen.Ma2024a}. Let $\mathbb{N}$ and $\mathbb{Z}^+$ be sets of all natural numbers and positive integers, respectively. The cardinal of a set $X$ is denoted by $|X|$. The power set of $X$ is denoted by $\mathcal{P}(X)$. We use Boolean operations $\cap$, $\cup$ and $\comb{\cdot}$ (complementation) on $\mathcal{P}(X)$.

\begin{definition}
    The {language of tense logic} consists of a denumerable set $\mathsf{Prop}=\{p_i: i\in\mathbb{N}\}$ of variables, connectives $\bot$ and $\to$, and unary tense operators $\B$ and $\bd$. The set $\L$ of all formulas is defined by:
    \[
    \L\ni \phi ::= p \mid \bot \mid (\phi\to\phi) \mid \B\phi\mid \bd\phi,~\text{where $p \in \mathsf{Prop}$.}
    \]
    The connectives $\top,\neg,\wedge$ and $\vee$ are defined as usual. Let $\D\phi:=\neg\B\neg\phi$ and $\bb\phi:=\neg\bd\neg\phi$. 
    The {\em modal degree} $md(\phi)$ of a formula $\phi$ is defined inductively as follows: 
    \begin{align*}
    md(p) &= 0 = md(\bot),\\
    md(\phi\to\psi)&=\max\{md(\phi),md(\psi)\},\\
    md(\B\phi)&=md(\phi)+1=md(\bd\phi).
    \end{align*}

    A {\em substitution} is a homomorphism $(.)^s:\L\twoheadrightarrow \L$ on the formula algebra $\L$.
\end{definition}

\begin{definition}
    A {\em frame} is a pair $\F =(X,R)$ where $X$ is a nonempty set and $R\sub X\times X$. We write $Rxy$ when $\tup{x,y}\in R$. The {\em inverse} of $R$ is defined as $\breve{R}=\{\tup{y,x}: Rxy\}$. For every $x\in X$, let $R[x]=\{y\in X: Rxy\}$ and $\breve{R}[x]=\{u\in X: Ryx\}$. For every $U\sub X$, we define $R[U]=\bigcup_{x\in U}R[x]$ and $\breve{R}[U]=\bigcup_{x\in U}\breve{R}[x]$. Let $\mathsf{Fr}$ and $\mathsf{Fin}$ denote the class of all frames and finite frames, respectively.

    A {\em general frame} is a triple $\gf=(X,R,A)$ where $(X,R)$ is a frame and $A\sub\mathcal{P}(X)$ is a set such that $\ve\in A$ and $A$ is closed under the operators $\cap$, $(\cdot)^c$, $R[\cdot]$ and $\breve{R}[\cdot]$. We call $A$ {\em the set of internal sets in $\gf$}. We write $\kappa\gf$ for the \emph{underlying frame} $(X,R)$ of $\gf$. Let $\mathsf{GF}$ denote the class of all general frames.
\end{definition}

\begin{definition}\label{def:tran}
    Let $\F=(X,R)$ be a frame. Then $R$ is (i) \emph{reflexive} if $\F\md\forall{x}(Rxx)$, (ii) \emph{symmetric} if $\F\md\forall{xy}(Rxy\to Ryx)$; and (iii) \emph{transitive} if $\F\md\forall{xyz}(Rxy\wedge Ryz\to Rxz)$. We call $R^r=R\cup\cset{\tup{x,x}:x\in X}$ the \emph{reflexive closure of $R$} and $R\cup\breve{R}$ the \emph{symmetric closure of $R$}. Note that $R^t=\bigcap\cset{R'\supseteq R:R'\text{ is transitive}}$ is the smallest transitive relation containing $R$, we call $R^t$ the \emph{transitive closure of $R$}. We say that $\F$ is reflexive and transitive if $R$ is reflexive and transitive, respectively. Let $\F^r=(X,R^r)$ and $\F^t=(X,R^t)$ be the reflexive and transitive closure of $\F$, respectively. 
\end{definition}

\begin{definition}
    Let $\F=(X,R)$ be a transitive frame. For each $x\in X$, the \emph{cluster generated by $x$} is defined as $C(x)=(R[x]\cap\breve{R}[x])\cup\cset{x}$. A subset $C\sub X$ is called a \emph{cluster in $\F$} if $C=C(x)$ for some $x\in X$. Then we say that $\F$ is a \emph{cluster} if $X$ is a cluster in $\F$. Moreover, $\F$ is a \emph{non-degenerated cluster} if $\F$ is a cluster and $R\neq\ve$.
\end{definition}

\begin{definition}
    Let $\gf=(X,R,A)$ be a general frame. Then a map $V:\mathsf{Prop}\to A$ is called a valuation in $\gf$. A valuation $V$ is extended to $V:\L\to A$ as follows:
    \begin{align*}
        V(\bot) &=\ve,
        &V(\phi\to\psi) &= (V(\phi))^c\cup V(\psi),\\
        V(\bd\phi) &= \bd_R V(\phi),
        &V(\B\phi) &= \B_R V(\phi),
    \end{align*}
    where operations $\B_R$ and $\bd_R$ are defined %
    by $\bd_R:U\mapsto R[U]$ and $\B_R:U\mapsto (\breve{R}[U^c])^c$. 
    
    A model is a pair $\M=(\gf,V)$ where $\gf\in\mathsf{GF}$ and $V$ a valuation in $\gf$. Let $\phi$ be a formula and $x\in X$. Then (i) $\phi$ is {\em true} at $x$ in $\M$ (notation: $\M,x\models\phi$) if $x\in V(\phi)$; (ii) $\phi$ is {\em valid} at $x$ in $\gf$ (notation: $\gf,x\models\phi$) if $x\in V(\phi)$ for every valuation $V$ in $\gf$; (iii) $\phi$ is {\em valid} in $\gf$ (notation: $\gf\models\phi$) if $\gf,x\models\phi$ for every $x\in X$; and (iv) $\phi$ is {\em valid} in a class $\mathcal{K}$ of general frames (notation: $\mathcal{K}\models\phi$) if $\gf\models\phi$ for every $\gf\in\mathcal{K}$. For each set $\Sigma\sub\L$ of formulas and class $\mathcal{K}\sub\mathsf{GF}$ of general frames, let 
    \begin{center}
        $\mathcal{K}(\Sigma)=\cset{\gf\in\mathcal{K}:\gf\md\Sigma}$ and $\mathsf{Log}(\mathcal{K})=\cset{\phi:\mathcal{K}\md\phi}$.
    \end{center}
\end{definition}

\begin{definition}
    Let $\gf=(X,R,A)$ be a general frame and $x\in X$. For $k\geq 0$, we define the set $R_\sharp^k[x]$ of $k$-reachable points from $x$ inductively as follows:
    \begin{center}
        $R_\sharp^0[x] = \{x\}$; $R_\sharp^{k+1}[x] =R_\sharp^k[x]\cup R[R_\sharp^k[x]]\cup \breve{R}[R_\sharp^k[x]]$.
    \end{center}
    Let $R_\sharp^\omega[x] = \bigcup_{k\geq 0}R_\sharp^k[x]$. An $\R$-path is a finite tuple $\tup{x_i:i\leq n}$ such that $x_i\in\R[x_{i+1}]$ for all $i<n$. A general frame $\gf=(X,R,A)$ is said to be \emph{rooted} or \emph{connected} if $X=\R^\omega[x]$ for some $x\in X$. 
\end{definition}

Clearly, a general frame $\gf=(X,R,A)$ is rooted if and only if $X=\R^\omega[x]$ for each $x\in X$.

\begin{definition}
    For each $n\in\omega$ and $\phi,\psi\in\L$, we define the formula $\Delta_\psi^n\phi$ by:
\begin{center}
    $\Delta_\psi^0\phi=\psi\wedge\phi$ and $\Delta_\psi^{k+1}\phi=\Delta_\psi^k\phi\vee\D(\psi\wedge\Delta_\psi^k\phi)\vee\bd(\psi\wedge\Delta_\psi^k\phi)$.
\end{center}
As usual, we define the dual operator $\nabla_\psi^n$ of $\Delta_\psi^n$ by $\nabla_\psi^n\phi:=\neg\Delta_\psi^n\neg\phi$. If $\psi=\top$, then we write $\Delta^n\phi$ and $\nabla^n\phi$ for $\Delta_\psi^n\phi$ and $\nabla_\psi^n\phi$, respectively.
\end{definition}

\begin{proposition}
    Let $\M=(X,R,V)$ be a model, $x\in X$ and $\phi,\psi\in\L$. Then for all $k\in\omega$,
    \begin{enumerate}[(1)]
        \item $\M,x\md\Delta_\psi^k\phi$ if and only if there exists an $\R$-path $\tup{x_i:i<k}$ such that $\M,x_{k-1}\md\phi$, $x=x_0$ and $\M,x_i\md\psi$ for all $i<k$.
        \item $\M,x\md\Delta^k\phi$ if and only if $\M,y\md\phi$ for some $y\in\R^k[x]$.
    \end{enumerate}
\end{proposition}
\begin{proof}
    By induction on $k$.
\end{proof}

\begin{definition}
    Let $\gf=(X,R,A)$ be a general frame. For every subset $Y$ of $X$, the {\em subframe of $\gf$ induced by $Y$} is defined as $\gf\rsto Y=(Y,R\rsto Y,A\rsto Y)$, where $R\rsto Y=R\cap(Y\times Y)$ and $A\rsto Y=\cset{U\cap Y:U\in A}$. Let $\gg$ be a general frame. If $\gg\cong\gf\rsto Y$, then we say $\gg$ can be embedded into $\gf$ and write $\gg\rightarrowtail\gf$. For all $x\in X$, let $\gf_x$ denote the frame $\gf\rsto \R^\omega[x]$. For each class $\mathcal{K}$ of general frames, let $\mathcal{K}_r=\cset{\gf_x:\gf\in\mathcal{K}\text{ and }x\in\gf}$.
    
    Let $\gf=(X,R,A)$ and $\gf'=(X',R',A')$ be general frames. A map $f:X\to X'$ is said to be a {\em t-morphism from $\gf$ to $\gf'$}, if $f^{-1}[Y']\in A$ for any $Y'\in A'$ and
    \begin{center}
        for all $x\in X$, $f[R[x]]= R'[f(x)]$ and $f[\breve{R}[x]]=\breve{R'}[f(x)]$.
    \end{center}
    We write $f:\gf\twoheadrightarrow\gf'$ ($f:\gf\cong\gf'$) if $f$ is a surjective (bijective) t-morphism from $\gf$ to $\gf'$. Moreover, $\gf'$ is called a {\em t-morphic (isomorphic) image of $\gf$} and we write $\gf\twoheadrightarrow\gf'$ ($\gf\cong\gf'$) if there exists $f:\gf\twoheadrightarrow\gf'$ ($f:\gf\cong\gf'$). For each class $\mathcal{K}$ of general frames, let $\mathsf{TM}(\mathcal{K})$ denote the class of all t-morphic images of frames in $\mathcal{K}$.
\end{definition}

\begin{remark}
    In this work, we identify isomorphic general frames, i.e., for all $\gf,\gf'\in\mathsf{GF}$, we say $\gf=\gf'$ if $\gf\iso\gf'$. Thus we are always allowed to rename the elements in our domain.
\end{remark}

\begin{proposition}
    Let $\gf=(X,R,A)\in\mathsf{GF}$, $x\in X$ and $\phi\in\L$. Then 
    \begin{center}
        $\gf,x\md\phi$ if and only if $\gf\rsto\R^{\mathsf{md}(\phi)},x\md\phi$.
    \end{center}
\end{proposition}
\begin{proof}
    By induction on $\mathsf{md}(\phi)$.
\end{proof}

\begin{definition}
    A {\em tense logic} is a set of formulas $L$ such that 
    (i) the classical propositional logic is a subset of $L$;
    (ii) $\bd \phi\to \psi\in L$ if and only if $\phi\to \B\psi\in L$;
    (iii) if $\phi,\phi\to\psi\in L$, then $\psi\in L$; 
    (iv) if $\phi\in L$, then $\phi^s\in L$ for all substitution $s$. 
    The least tense logic is denoted by $\mathsf{K}_t$. 
\end{definition}
    
For every tense logic $L$ and set of formulas $\Sigma$, let $L\oplus\Sigma$ denote the smallest tense logic containing $L\cup\Sigma$. A tense logic $L_1$ is a {\em sublogic} of $L_2$ (or $L_2$ is an {\em extension} of $L_1$) if $L_1\sub L_2$. Let $\mathsf{NExt}(L)$ be the set of all extensions of $L$. The readers can check that for all tense logic $L$, $(\mathsf{NExt}(L), \cap,\oplus)$ is a distributive lattice with top $\L$ and bottom $L$. Moreover, since tense logics are closed under arbitrary intersections, $\mathsf{NExt}(L)$ is a complete lattice. For each subset $X\sub\mathsf{NExt}(L)$, we write $\bigoplus X$ for the supreme of $X$. The following completeness result is well-known:

\begin{theorem}\label{thm:gf-completeness}
    Let $L$ be a tense logic. Then $L=\mathsf{Log}(\mathsf{GF}(L))=\mathsf{Log}(\mathsf{GF}_r(L))$.
\end{theorem}

\begin{lemma}\label{lem:intersection-rooted-frame}
    Let $L_1,L_2$ be tense logics and $\mathcal{K}\sub\mathsf{GF}_r$. Then $\mathcal{K}(L_1\cap L_2)=\mathcal{K}(L_1)\cup\mathcal{K}(L_2)$.
\end{lemma}
\begin{proof}
    Clearly $\mathcal{K}(L_1\cap L_2)\supseteq\mathcal{K}(L_1)\cup\mathcal{K}(L_2)$. Take any $\gf=(X,R,A)\in\mathcal{K}(L_1\cap L_2)$. Suppose $\gf\not\in\mathcal{K}(L_1)\cup\mathcal{K}(L_2)$. Then there are $\psi_1\in L_1$, $\psi_2\in L_2$ and $x,y\in X$ such that $\gf,x\not\md\psi_1$ and $\gf,y\not\md\psi_2$. Moreover, we may assume there is no common variable in $\psi_1$ and $\psi_2$. Thus there exists a valuation $V$ in $\gf$ such that $\gf,V,x\md\neg\psi_1$ and $\gf,V,y\md\neg\psi_2$. Since $\gf$ is rooted, $y\in\R^n[x]$ for some $n\in\omega$. Thus $\gf,V,x\md\neg\psi_1\wedge\neg\nabla^k\psi_2$ Note that $\psi_1\vee\nabla^k\psi_2\in L_1\cap L_2$, we have $\gf\not\in\mathcal{K}(L_1\cap L_2)$, which leads to a contradiction.
\end{proof}

\begin{definition}
    Let $L$ be a tense logic. Then we say (i) $L$ is \emph{Kripke complete}, if $L=\mathsf{Log}(\mathsf{Fr}(L))$; (ii) $L$ enjoys the \emph{finite model property (FMP)}, if $L=\mathsf{Log}(\mathsf{Fin}(L))$.
\end{definition}

\begin{definition}
    Let $(\mathsf{T})$, $(\mathsf{4})$, $(\mathsf{5})$, $(\mathsf{grz}^+)$ and $(\mathsf{grz}^-)$ denote the following formulas:
    \begin{align*}
        \tag*{($\mathsf{T}$)} \B p\to p\\
        \tag*{($\mathsf{4}$)} \B p\to\B\B p\\
        \tag{$\mathsf{5}$} \D p\to\B\D p
    \end{align*}
    Let $\mathsf{K4}_t=\mathsf{K}_t\oplus\mathsf{4}$, $\ST=\mathsf{K4}_t\oplus\mathsf{T}$ and $\mathsf{S5}_t=\ST\oplus\mathsf{5}$.
\end{definition}

\begin{fact}
    Let $\gf=(X,R,A)\in\mathsf{RF}_r$. Then (i) $\gf\md\mathsf{K4}_t$ if and only if $\gf$ is transitive; (ii) $\gf\md\mathsf{S4}_t$ if and only if $\gf$ is transitive and reflexive; (iii) $\gf\md\mathsf{S5}_t$ if and only if $\gf$ is a non-degenerated cluster. 
\end{fact}

\begin{definition}
    Let $\gf=(X,R,A)\in\mathsf{GF}$ be a frame and $x\in X$. Let $k\in\mathbb{Z}^+$. Then we say {\em $x$ is of reachability-degree (r-degree) $k$ (notation: $\mathrm{rdg}(x)=k$)}, if $\R^{k-1}[x]\neq \R^{k}[x]=\R^\omega[x]$. Specially, $\mathrm{rdg}(x)=0$ if $\R^\omega[x]=\cset{x}$ and $\mathrm{rdg}(x)=\aleph_0$ if $\R^k[x]\neq\R^{k+1}[x]$ for any $k\in\omega$. We define the {\em r-degree $\mathrm{rdg}(\gf)$ of $\gf$} by $\mathrm{rdg}(\gf)=\mathrm{sup}\cset{\mathrm{rdg}(x):x\in X}$. Then we say $\gf$ is \emph{$n$-transitive} if $\mathsf{rdg}(\gf)\leq n$. A tense logic $L$ is (i) $n$-transitive if each of its general frames is $n$-transitive; (ii) finitely transitive if it is $n$-transitive for some $n\in\omega$. 
\end{definition}

\begin{remark}
    A normal modal logic $L$ is called pre-transitive if $\B^np\to\B^{n+1}p\in L$ for some $n\in\omega$. The notion finitely transitivity in the tense case is similar to pre-transitivity in the modal case, which guarantee the existence of a master modality, say $\Delta^n$. A tense logic $L$ is $n$-transitive if and only if $\Delta^{n+1}p\to\Delta^n p\in L$. It is worth to notice that the tense logic $\mathsf{S4}_t$ is not finitely transitive.
\end{remark}

\begin{definition}
    Let $\gf=(X,R,A)\in\mathsf{GF}$, $\alpha$ an ordinal and $\mathcal{Y}=\tup{y_i\in X:i<\alpha}$. Then 
    \begin{itemize}
        \item $\mathcal{Y}$ is called a chain in $\F$ if $Rx_\lambda x_\gamma$ for all $\lambda<\gamma<\alpha$;
        \item $\mathcal{Y}$ is called a strict chain in $\F$ if it is a chain and $x_\lambda\not\in R[x_\gamma]$ for all $\lambda<\gamma<\alpha$;
        \item $\mathcal{Y}$ is called a (strict) co-chain in $\F$ if it is a (strict) chain in $\breve{\F}$;
        \item $\cset{y_i\in X:i<\alpha}$ is called an anti-chain in $\F$ if $x_\lambda\not\in R[x_\gamma]$ for all $\lambda\neq\gamma<\alpha$.
    \end{itemize}
    The length $l(\mathcal{Y})$ of a strict chain $\mathcal{Y}=\tup{y_i\in X:i<\alpha}$ is defined to be $\alpha$.
    We say that $x\in X$ is of depth $n$ (notation: $\mathrm{dep}(x)=n$), if there exists a strict chain $\mathcal{Y}$ in $\gf\rsto R[x]$ with $l(\mathcal{Y})=n$ and there is no strict chain of greater length. Otherwise $x$ is said to be of infinite depth and we write $\mathrm{dep}(x)=\aleph_0$. We define the {\em depth $\mathsf{dep}(\gf)$ of $\gf$} by $\mathsf{dep}(\gf)=\mathrm{sup}\cset{\mathsf{dep}(x):x\in X}$.

    Let $n\in\mathbb{Z}^+$. We say that $x\in X$ is of forth-width $n$ (notation: $\mathrm{wid}^+(x)=n$), if there exists an anti-chain $Y\sub R[x]$ with $|Y|=n$ and there is no anti-chain in $R[x]$ with greater size. Otherwise we write $\mathrm{wid}^+(\gf)=\aleph_0$. Back-width is defined dually and we write $\mathrm{wid}^-(x)=n$ if $x$ is of back-width $n$.
\end{definition}

\begin{definition}
    For each $n\in\mathbb{Z}^+$, we define the formulas $(\mathsf{alt}^+_n)$, $(\mathsf{alt}^-_n)$, $(\mathsf{bz}_n)$, $(\mathsf{bw}^+_n)$ and $(\mathsf{bw}^-_n)$ as follows:
    \begin{align*}
        \tag*{($\mathsf{alt}^+_n$)} \B p_0 \vee \B(p_0 \to p_1) \vee \cdots \vee \B(p_0 \wedge \cdots \wedge p_{n-1} \to p_{n})\\
        \tag*{($\mathsf{alt}^-_n$)} \bb p_0 \vee \bb(p_0 \to p_1) \vee \cdots \vee \bb(p_0 \wedge \cdots \wedge p_{n-1} \to p_{n})\\
        \tag{$\mathsf{bz}_n$} \Delta^{n+1}p\to\Delta^{n}p\\
        \tag{$\mathsf{bw}^+_n$} \bigwedge_{i\leq n}\D p_i\to\bigvee_{i\neq j\leq n}\D(p_i\wedge(p_j\vee\D p_j))\\
        \tag{$\mathsf{bw}^-_n$} \bigwedge_{i\leq n}\bd p_i\to\bigvee_{i\neq j\leq n}\bd(p_i\wedge(p_j\vee\bd p_j))
    \end{align*}
    Moreover, we define the formula $(\mathsf{bd}_n)$ for each $n\in\mathbb{Z}^+$ as follow:
    \begin{align*}
        \mathsf{bd}_1 &= \D\B p_0\to p_0\\
        \mathsf{bd}_{k+1} &= \D(\B p_{k}\wedge\neg\mathsf{bd}_k)\to p_{k}
    \end{align*}
\end{definition}

\begin{fact}\label{fact:bounds}
    Let $\gf=(X,R,A)\in\mathsf{RF}$, $x\in X$ and $n\in\mathbb{Z}^+$. Then
    \begin{enumerate}[(1)]
        \item $\gf,x\md\mathsf{alt}^+_n$ if and only if $|R[x]|\leq n$.
        \item $\gf,x\md\mathsf{alt}^-_n$ if and only if $|\breve{R}[x]|\leq n$.
        \item $\gf,x\md\mathsf{bz}_n$ if and only if $\mathrm{rdg}(x)\leq n$.
        \item $\gf,x\md\mathsf{bw}^+_n$ if and only if $\mathrm{wid}^+(x)\leq n$.
        \item $\gf,x\md\mathsf{bw}^-_n$ if and only if $\mathrm{wid}^-(x)\leq n$.
        \item $\gf,x\md\mathsf{bd}_n$ if and only if $\mathrm{dep}(x)\leq n$.
    \end{enumerate}
\end{fact}
\begin{proof}
    The proof for this fact is standard, see \cite{Chagrov.Zakharyaschev1997}.
\end{proof}

\begin{definition}
    Let $(\mathsf{grz}^+)$ and $(\mathsf{grz}^-)$ denote the formulas $\B(\B(p\to\B p)\to p)\to p$ and $\bb(\bb(p\to\bb p)\to p)\to p$, respectively. %
\end{definition}

\begin{fact}\label{fact:grz}
    Let $\F=(X,R)\in\mathsf{Fr}(\ST)$ and $x\in X$. Then (i) $\F,x\md\mathsf{grz}^+$ if and only if $\breve{R}$ is well-founded, i.e., there is no infinite ascending chain in $\F$; (ii) $\F,x\md\mathsf{grz}^-$ if and only if $\breve{\F},x\md\mathsf{grz}^+$.
\end{fact}

\begin{definition}
    Let $\F=(X,R)$ be a finite rooted frame and $x\in X$. Let $k\in\mathbb{Z}^+$ and $\tup{x_i:i\in n}$ be an enumeration of $X$. The formula $\J^k(\F)$ is defined to be the conjunction of the following formulas:
    \begin{enumerate}[(1)]
        \item $\nabla^k(p_0\vee\cdots\vee p_{n-1})$
        \item $\nabla^k(p_i\to\neg p_j)$, for all $i\neq j$
        \item $\nabla^{k-1}((p_i\to\D p_j)\wedge(p_j\to\bd p_i))$, for all $\tup{x_i,x_j}\in R$
        \item $\nabla^{k-1}((p_i\to\neg\D p_j)\wedge(p_j\to\neg\bd p_i))$, for all $\tup{x_i,x_j}\not\in R$
    \end{enumerate}
    $\J^k(\F)$ is called the Jankov formula of $\F$ of degree $k$.
\end{definition}

It is not hard to check that $\F\not\md\neg\J^k(\F)$ for any $k\in\mathbb{Z}^+$. Similar to the case for pre-transitive modal logic, the following lemma holds:

\begin{theorem}[{\cite[Theorem 3.12]{Chen2024}}]\label{thm:JankovLemma}
    Let $k\in\mathbb{Z}^+$, $\gf=(X,R,A)$ be a $k$-transitive general frame and $\G=(Y,S)$ a finite rooted frame. Then
    \begin{center}
        $\F\not\md\neg\J^k(\G)$ if and only if $\gf_x\twoheadrightarrow\G$ for some $x\in X$.
    \end{center}
\end{theorem}

\begin{definition}
    Let $L_0$ be a tense logic and $L\in\mathsf{NExt}(L_0)$. Then we define the degree of Kripke-incompleteness $\mathsf{deg}_{\mathsf{NExt}(L_0)}(L)$ of $L$ in $\mathsf{NExt}(L_0)$ by
\begin{center}
    $\mathsf{deg}_{\mathsf{NExt}(L_0)}(L)=|\cset{L'\in\mathsf{NExt}(L_0):\mathsf{Fr}(L')=\mathsf{Fr}(L)}|$.
\end{center}
To simplify notation, we write $\mathsf{deg}_{L_0}$ for $\mathsf{deg}_{\mathsf{NExt}(L_0)}$. Moreover, we define the degree of finite model property (degree of FMP) $\mathsf{df}_{\mathsf{NExt}(L_0)}(L)$ of $L$ in $\mathsf{NExt}(L_0)$ by
\begin{center}
    $\mathsf{df}_{\mathsf{NExt}(L_0)}(L)=|\cset{L'\in\mathsf{NExt}(L_0):\mathsf{Fin}(L')=\mathsf{Fin}(L)}|$.
\end{center}
Again, we write $\mathsf{df}_{L_0}$ for $\mathsf{df}_{\mathsf{NExt}(L_0)}$.
\end{definition}

The following proposition follows immediately from the fact that $\mathsf{Fin}(L)\sub\mathsf{Fr}(L)$:
\begin{proposition}\label{prop:deg-leq-df}
    Let $L_0$ be a tense logic and $L\in\mathsf{NExt}(L_0)$. Then $\mathsf{deg}_{L_0}(L)\leq\mathsf{df}_{L_0}(L)$.
\end{proposition}

\begin{definition}
    Let $L_0\in\NExt(\mathsf{K}_t)$ and $L_1,L_2\in\mathsf{NExt}(L_0)$. Then $\tup{L_1,L_2}$ is called a splitting pair in $\mathsf{NExt}(L_0)$ if, for all $L\in\mathsf{NExt}(L_0)$, exactly one of $L\sub L_1$ and $L\supseteq L_2$ holds. In this case, we say $L_1$ splits the lattice $\mathsf{NExt}(L_0)$ and we call $L_2$ the splitting of $\mathsf{NExt}(L_0)$ by $L_1$ and denote it by $L_0/L_1$. 
    
    A tense logic $L$ is called a union splitting in $\mathsf{NExt}(L_0)$ if there exists a family $\cset{L_i:i\in I}$ of splittings in $\mathsf{NExt}(L_0)$ such that $L=\bigoplus_{i\in I}L_i$. We say that $L$ is an iterated splitting in $\NExt(L_0)$ if $L=L_0/L_1/\cdots/L_n$ for some $L_1,\cdots,L_n$ such that for all $1\leq i\leq n$, $L_i$ splits $\NExt(L_0/L_1/\cdots/L_{i-1})$. Specially, we count also $L_0$ itself an iterated splitting of $\NExt(L_0)$.
\end{definition}

Splittings of lattices of logics play a core role in Blok's characterization of degree of Kripke-incompleteness. Consider the lattice $\mathsf{NExt}(\mathsf{K})$ of normal modal logics. It is known that a normal modal logic $L$ splits $\mathsf{NExt}(\mathsf{K})$ if and only if $L=\mathsf{Log}(\F)$ for some finite rooted cycle free frame $\F$ (see \cite[Theorems 10.49 and 10.53]{Chagrov.Zakharyaschev1997}). Blok \cite{Blok1978a} proved that a normal modal logic $L$ has the degree of incompleteness $1$ in $\mathsf{NExt}(\mathsf{K})$ if and only if $L$ is a union splitting of $\mathsf{NExt}(\mathsf{K})$. It was also showed in \cite{Blok1978a} that in the lattice $\NExt(\mathsf{K})$, consistent iterated splittings are exactly union-splittings.

Let us first recall some general results on splitting pairs in lattices from general lattice theory.

\begin{fact}\label{fact:splitting}
    Let $\mathcal{L}=\tup{L,\leq}$ be a lattice and $x\in L$. Then 
    \begin{enumerate}[(1)]
        \item $x$ splits $\mathcal{L}$ if and only if $x$ is completely $\wedge$-prime in $\mathcal{L}$, i.e., for all $\cset{y_i:i\in I}\sub L$, $\bigwedge_{i\in I}y_i\leq x$ implies $y_i\leq x$ for some $i\in I$.
        \item $x$ is a splitting in $\mathcal{L}$ if and only if $x$ is completely $\vee$-prime in $\mathcal{L}$, i.e., for all $\cset{y_i:i\in I}\sub L$, $x\leq\bigvee_{i\in I}y_i$ implies $x\leq y_i$ for some $i\in I$.
    \end{enumerate}
\end{fact}

Lattices of tense logics are substantially different from those of modal logics. Splittings of lattices of tense logics have been investigated in \cite{Kracht1992}. Let $\bullet$ denote the frame $(\cset{0},\ve)$. Then

\begin{theorem}[{\cite[Theorem 22]{Kracht1992}}]\label{thm:splittings-Kt-K4t}
    Let $L\in\mathsf{NExt}(\mathsf{K}_t)$. Then
    \begin{enumerate}[(1)]
        \item $L$ splits $\mathsf{NExt}(\mathsf{K}_t)$ if and only if $L=\mathsf{Log}(\bullet)$.
        \item $L$ splits $\mathsf{NExt}(\mathsf{K4}_t)$ if and only if $L=\mathsf{Log}(\bullet)$.
    \end{enumerate}
\end{theorem}
It turns out that there exists exactly one splitting pair in $\mathsf{NExt}(\mathsf{K}_t)$, while there are countably many in $\mathsf{NExt}(\mathsf{K})$. However, as we are going to prove in Section~\ref{sec:DKI-Kt-K4t}, union-splittings in $\K$ and $\LT$ are still exactly the strictly Kripke-complete logics. Moreover, as Blok \cite{Blok1978a} proved that consistent iterated splittings in $\NExt(\mathsf{K})$ are exactly the union-splittings, we will show in Section~4 that this holds for also $\K$ and $\LT$. However, in Section~\ref{sec:DKI-S4t}, we show that there exists iterated splittings which are not union-splittings in $\NExt(\ST)$.

\section{Reflective Unfolding of Kripke frames}

In this section, we introduce the reflective unfolding of Kripke frames, which was studied by Kracht \cite{Kracht1992}. The reflective unfolding turns out to be one of the crucial constructions in our proofs of the characterization theorems. In our proofs, to construct a set of pairwise different tense logics which shares the same class of frames, as it will be explained in the following sections, we need frames of large enough r-degree. In general, we could obtain such frame by the unrevealing method (see \cite{Blackburn.deRijke.ea2001}). However, if we require further that these frames are transitive, then unrevealing does not always suffice. This is precisely why we need the reflective unfolding.

\begin{definition}
    Let $\F=(X,R)$ and $\G=(Y,S)$ be frames, $w\in X$ and $u\in Y$. Then we define the \emph{combination $\tup{\F{w{+}u}\G}=(Z,T)$ of $\F$ and $\G$ at $\tup{w,u}$} by $Z=(X\times\cset{0})\cup(Y\times\cset{1})\setminus\cset{\tup{u,1}}$ and
    \begin{align*}
        T=&\cset{\tup{\tup{x,a},\tup{y,a}}\in Z\times Z:\tup{x,y}\in R\cup S\text{ and }a\in\cset{0,1}}\cup\\
        &(\cset{\tup{w,0}}\times\cset{\tup{x,1}:x\in S[u]})\cup(\cset{\tup{x,1}:x\in\breve{S}[u]}\times\cset{\tup{w,0}})
    \end{align*}
    We call $\tup{\F{w{+^t}u}\G}=(Z,T^t)$ the transitive combination of $\F$ and $\G$ at $\tup{w,u}$. 
\end{definition}

\begin{example}
    Let $\F=(X,R)$ and $\G=(Y,S)$ be frames, $w\in X$ and $u\in Y$. The first observation is that if $w=u$ and $X\cap Y=\cset{u}$, then $\tup{\F{u{+}u}\G}=(X\cup Y,R\cup S)$. Intuitively, the frame $\tup{\F{w{+}u}\G}=(Z,T)$ is obtained by taking the disjoint union of $\F$ and $\G$ and then identifying the points $w$ and $u$. To see this, we give another equivalent definition of $T$ as follows
    \begin{itemize}
        \item for all $x,y\in X$, $T\tup{x,0}\tup{y,0}$ if and only if $Rxy$;
        \item for all $x,y\in Y\setminus\cset{u}$, $T\tup{x,1}\tup{y,1}$ if and only if $Sxy$;
        \item for all $y\in Y\setminus\cset{u}$, $T\tup{w,0}\tup{y,1}$ if and only if $Suy$;
        \item for all $y\in Y\setminus\cset{u}$, $T\tup{y,1}\tup{w,0}$ if and only if $Syu$.
    \end{itemize}
    Recall that we identify isomorphic frames, it is always possible to assume that $X\cap Y=\ve$. In this case, the definition of $\tup{\F{w{+}u}\G}=(Z,T)$ can be simplified as follows: $Z=(X\cup Y)\setminus\cset{u}$ and $T=(R\cup S\cup(\cset{w}\times S[u])\cup(\breve{S}[u]\times\cset{w}))\cap(Z\times Z)$. In what follows, we always presume that domains of frames are disjoint and use the simplified definition.    
\end{example}

\begin{fact}\label{fact:plus-reachability}
    Let $\F=(X,R)$ and $\G=(Y,S)$ be frames, $w\in X$ and $u\in Y$. Let $\tup{\F{w{+^t}u}\G}=(Z,T^t)$. Then for all $x,x'\in X$ and $y,y'\in Y\setminus\cset{u}$, the following holds:
    \begin{enumerate}[(1)]
        \item $R^txx'$ if and only if $T^txx'$;
        \item $S^tyy'$ if and only if $T^tyy'$;
        \item $T^txy$, if and only if, $(R^t)^rxw$ and $S^tuy$;
        \item $\breve{T^t}xy$, if and only if, $(\breve{R}^t)^rxw$ and $\breve{S}^tuy$.
    \end{enumerate}
\end{fact}
\begin{proof}
    Follows from the construction $\tup{\F{w{+}u}\G}$ immediately.
\end{proof}

\begin{figure}[ht]
\[
    \begin{tikzpicture}[scale=0.6]
        \draw (1,-3) node{$\F$};
        \draw (0,0) node{$\bullet$};
        \draw (0,0) node[above]{\small$w$};
        \draw [->] (.05,-.05) -- (0.95,-.95);
        \draw (1,-1) node{$\bullet$};
        \draw (1,-1) node[below]{\small$v$};
        \draw [->] (.05,.05) -- (1.95,.95);
        \draw (2,1) node{$\bullet$};
        \draw (2,1) node[above]{\small$u$};
    \end{tikzpicture}
    \quad
    \begin{tikzpicture}[scale=0.6]
        \draw (1,-3) node{$\F^1_{w,u}$};
        \draw (0,0) node{$\bullet$};
        \draw (0,0) node[above]{\small$w_0$};
        \draw [->] (.05,-.05) -- (0.95,-.95);
        \draw (1,-1) node{$\bullet$};
        \draw (1,-1) node[below]{\small$v_0$};
        \draw [->] (.05,.05) -- (1.95,.95);
        \draw (2,1) node{$\bullet$};
        \draw (2,1) node[above]{\small$u_0$};
    \end{tikzpicture}
    \quad
    \begin{tikzpicture}[scale=0.6]
        \draw (2,-3) node{$\F^2_{w,u}$};
        \draw (0,0) node{$\bullet$};
        \draw (0,0) node[above]{\small$w_0$};
        \draw [->] (.05,-.05) -- (0.95,-.95);
        \draw (1,-1) node{$\bullet$};
        \draw (1,-1) node[below]{\small$v_0$};
        \draw [->] (.05,.05) -- (1.95,.95);
        \draw (2,1) node{$\bullet$};
        \draw (2,1) node[above]{\small$u_0$};
        \draw [<-] (2.05,.95) -- (3.95,.05);
        \draw (4,0) node{$\bullet$};
        \draw (4,0) node[above]{\small$w_1$};
        \draw [->] (3.95,-.05) -- (3.05,-.95);
        \draw (3,-1) node{$\bullet$};
        \draw (3,-1) node[below]{\small$v_1$};
    \end{tikzpicture}
    \quad
    \begin{tikzpicture}[scale=0.6]
        \draw (4,-3) node{$\F^4_{w,u}$};
        \draw (0,0) node{$\bullet$};
        \draw (0,0) node[above]{\small$w_0$};
        \draw [->] (.05,-.05) -- (0.95,-.95);
        \draw (1,-1) node{$\bullet$};
        \draw (.5,-1) node[below]{\small$v_0$};
        \draw [->] (.05,.05) -- (1.95,.95);
        \draw (2,1) node{$\bullet$};
        \draw (2,1) node[above]{\small$u_0$};
        \draw [<-] (2.05,.95) -- (3.95,.05);
        \draw (4,0) node{$\bullet$};
        \draw (4,0) node[above]{\small$w_1$};
        \draw [->] (3.95,-.05) -- (3.05,-.95);
        \draw (3,-1) node{$\bullet$};
        \draw (3,-1) node[below]{\small$v_1$};
        \draw [->] (.05+4,-.05) -- (0.95+4,-.95);
        \draw (1+4,-1) node{$\bullet$};
        \draw (1+4,-1) node[below]{\small$v_2$};
        \draw [->] (.05+4,.05) -- (1.95+4,.95);
        \draw (2+4,1) node{$\bullet$};
        \draw (2+4,1) node[above]{\small$u_2$};
        \draw [<-] (2.05+4,.95) -- (3.95+4,.05);
        \draw (4+4,0) node{$\bullet$};
        \draw (4+4,0) node[above]{\small$w_3$};
        \draw [->] (3.95+4,-.05) -- (3.05+4,-.95);
        \draw (3+4,-1) node{$\bullet$};
        \draw (3+4,-1) node[below]{\small$v_3$};
    \end{tikzpicture}
\]
\caption{Examples for the reflective unfolding}
\label{fig:book}
\end{figure}
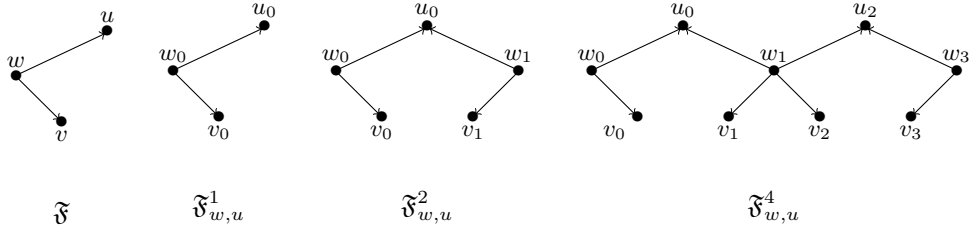

\begin{definition}
    Let $\F=(X,R)$ be a frame. For each $n\in\omega$, we write $\F[n]$ for the frame $(X[n],R[n])$, where $X[n]=\cset{x_n:x\in X}$ and $R[n]=\cset{\tup{x_n,y_n}:xRy}$. For all $w,u\in X$ and $n\in\mathbb{Z}^+$, we define the $n$-r-unfolding $\F^n_{w,u}=(X^n_{w,u},R^n_{w,u})$ of $\F$ by $(w,u)$ inductively as follows:
    \begin{itemize}
        \item $\F^1_{w,u}=\F[0]$;
        \item $\F^{2k+2}_{w,u}=\tup{\F^{2k+1}_{w,u}u_{2k}+u_{2k+1}\F[2k+1]}$;
        \item $\F^{2k+3}_{w,u}=\tup{\F^{2k}_{w,u}w_{2k+1}+w_{2k+2}\F[2k+2]}$.
    \end{itemize}
\end{definition}
Intuitively, the frame $\F^n_{w,u}$ is constructed by first get $n$ copies of $\F$ and then combine them in a special way. It should be clear that the reflective unfolding preserves reflexivity and connectedness, i.e., $\F^n_{w,u}$ is rooted and reflexive if $\F$ is rooted and reflexive, respectively. An example of the reflective unfolding is given in Figure \ref{fig:book}.

\begin{lemma}\label{lem:book-morphism}
    Let $\F=(X,R)$ be a frame and $w,u\in X$. For each $n\in\mathbb{Z}^+$, let $n\phi:X^n_{w,u}\to X$ be the map defined by $n\phi:x_k\mapsto x$ for all $x\in X$ and $k<n$. Then $n\phi$ is a t-morphism from $\F^n_{w,u}$ to $\F$. Moreover, if $\F$ is transitive, then $n\phi:(\F^n_{w,u})^t\twoheadrightarrow\F$.
\end{lemma}
\begin{proof}
    For all $x_i\in X^n_{w,u}$, $n\phi[R^n_{w,u}[x_i]]=R[n\phi(x_i)]=R[x]$ and $n\phi[\breve{R}^n_{w,u}[x_i]]=\breve{R}[n\phi(x_i)]=\breve{R}[x]$.
\end{proof}

\begin{lemma}\label{lem:book-degree}
    Let $\F=(X,R)\in\mathsf{Fr}_r$, $x,w,u\in X$ and $k,n\in\omega$ such that $w\neq u$ and $k<4n+2$. Let $\F^{4n+2}_{w,u}=(Y,S)$. Then the following holds:
    \begin{enumerate}[(1)]
        \item $S_\sharp[x_k]\sub\bigcup\cset{X[i]:i<4n+2,k\leq i+1\text{ and }i\leq k+1}$. 
        \item $S_\sharp^n[x_k]\neq Y$.
        \item if $\tup{w,u}\not\in R^t$, then $(S^t)_\sharp[x_k]\sub\bigcup\cset{X[i]:i<4n+2,k\leq i+2\text{ and }i\leq k+2}$. 
        \item if $\tup{w,u}\not\in R^t$, then $(S^t)_\sharp^n[x_k]\neq Y$.
    \end{enumerate}
\end{lemma}
\begin{proof}
    For (1), suppose $k$ is even. If $k=0$, then $S_\sharp[x_k]\sub X[0]\cup X[1]$. If $k>0$, then $S_\sharp[x_k]\sub\cset{w_{k-1}}\cup X[k]\cup (R[k+1])_\sharp[u_{k+1}]\sub X[k-1]\cup X[k]\cup X[k+1]$. Suppose $k$ is odd. Then $S_\sharp[x_k]\sub\cset{u_{k-1}}\cup X[k]\cup (R[k+1])_\sharp[w_{k+1}]\sub X[k-1]\cup X[k]\cup X[k+1]$. Thus (1) holds.
    By (1), we see $S_\sharp^n[x_k]\sub\bigcup\cset{X[i]:i<4n+2,k\leq n+i\text{ and }i\leq k+n}\neq Y$, from which (2) follows immediately.

    For (3), suppose $k$ is odd. Take any $y_j\in (S^t)_\sharp[x_k]$. It suffices to show that $k\leq j+2$ and $j\leq k+2$. Suppose $j>k+2$. Since $y_j\in (S^t)_\sharp[x_k]$, either $S^tx_ky_j$ or $S^ty_jx_k$. Suppose $S^tx_ky_j$. By Fact~\ref{fact:plus-reachability}, we see $(R^t)^rxu$ and $S^tu_ky_j$. By applying Fact~\ref{fact:plus-reachability} twice, we have $(R^t)^ruw$ and $(R^t)^rwu$. Since $w\neq u$, we have $R^twu$, which contradicts $\tup{w,u}\not\in R^t$. By Fact~\ref{fact:plus-reachability}, $S^ty_jx_k$ implies $(\breve{R}^t)^ruw$ and so $R^twu$, which also gives a contradiction. Thus $j\leq k+2$. By a similar argument, we see that $k\leq j+2$. The proof for the case when $k$ is exactly the dual one. Thus (3) holds. By (3), we see $S_\sharp^n[x_k]\sub\bigcup\cset{X[i]:i<4n+2,k\leq 2n+i\text{ and }i\leq k+2n}\neq Y$, (4) follows immediately.
\end{proof}

\begin{corollary}\label{coro:Kt-large-rdg}
    Let $\phi\in\L$. If $\phi$ is satisfied by some finite rooted frame other than $\bullet$, then for each $n\in\omega$, there exists $\F=(X,R)\in\mathsf{Fin}_r(L)$ such that $\F\not\md\neg\phi$ and $\mathsf{rdg}(\F)\geq n$.
\end{corollary}
\begin{proof}
    Suppose $\phi$ is satisfied by some rooted frame $\G=(Y,S)\neq\bullet$. Suppose $|Y|=1$. Then $S=Y\times Y$. Let $\F$ be a finite $\ST$-frame of r-degree greater than $n$, for example, $\G_{n+1}$ defined in \cite[Definition 7.1]{Chen2024}. Note that $\F\twoheadrightarrow\G$ and $\G\not\md\neg\phi$, we see that $\F\not\md\neg\phi$. Suppose $|Y|\geq 2$. Then there exists $w,u,y\in Y$ such that $\G,y\not\md\neg\phi$ and $w\neq u$. By Lemmas~\ref{lem:book-morphism} and \ref{lem:book-degree}(2), we see that $\G^{4n+2}_{w,u},y_0\not\md\neg\phi$ and $\mathsf{rdg}(\G^{4n+2}_{w,u})\geq n$.
\end{proof}

\begin{corollary}\label{coro:Tr-large-rdg}
    Let $\phi\in\L$ and $L\in\cset{\mathsf{K4}_t,\mathsf{S4}_t}$. If $\phi$ is satisfied by some rooted non-symmetric $L$-frame, then for each $n\in\omega$, there exists $\F=(X,R)\in\mathsf{Fin}_r(L)$ such that $\F\not\md\neg\phi$ and $\mathsf{rdg}(\F)\geq n$.
\end{corollary}
\begin{proof}
    Let $\G=(Y,S)$ be a non-symmetric $L$-frame and $y\in Y$ such that $\G,y\not\md\neg\phi$. Then there exists $w\neq u$ with $\tup{w,u}\not\in R$. By Lemmas~\ref{lem:book-morphism} and \ref{lem:book-degree}(4), we see that $\G^{4n+2}_{w,u},y_0\not\md\neg\phi$ and $\mathsf{rdg}(\G^{4n+2}_{w,u})\geq n$. Note that $\G^{4n+2}_{w,u}\md L$ whenever $\G\md L$, take $\F=\G^{4n+2}_{w,u}$ and we are done.
\end{proof}

\section{Kripke-incompleteness in $\K$ and $\LT$}\label{sec:DKI-Kt-K4t}

In this section, we generalize Blok's dichotomy theorem for $\mathsf{NExt}(\mathsf{K})$ to the lattices $\mathsf{NExt}(\mathsf{K}_t)$ and $\LT$ of tense logics. For $\mathsf{NExt}(\mathsf{K}_t)$ and $\LT$, we prove the dichotomy theorems by showing that the union-splittings of $\K$ and $\LT$ are exactly those having the degree of Kripke-incompleteness $1$, and all other tense logics have the degree of Kripke-incompleteness $2^{\aleph_0}$, respectively.

\subsection{Degree of Kripke-incompleteness in $\K$}\label{subsec:degree-Kt}
Let us focus on the degree of Kripke-incompleteness in $\K$ and we write $\mathsf{deg}$ for $\mathsf{deg}_{\mathsf{K}_t}$ in this section. By Theorem~\ref{thm:splittings-Kt-K4t}, $\tup{\mathsf{Log}(\bullet),\mathsf{K}_t/\mathsf{Log}(\bullet)}$ is the unique splitting pair in $\K$. Clearly, $\mathsf{Log}(\bullet)=\mathsf{K}_t\oplus(\B\bot\wedge\bb\bot)$. To simplify our notation, in this section, we write $L^*$ for $\mathsf{K}_t/\mathsf{Log}(\bullet)$. The reader can readily check that $\cset{\mathsf{K}_t,L^*}$ is exactly the set of union-splittings of $\K$. Moreover, since $\mathsf{K}_t$ is Kripke-complete, we see immediately that $\mathsf{deg}(\mathsf{K}_t)=1$.

\begin{proposition}\label{prop:deadend-frame}
    Let $\gf\in\mathsf{GF}_r$. Then $\gf\md\D\top\vee\bd\top$ if and only if $\gf\ncong\bullet$.
\end{proposition}
\begin{proof}
    The left-to-right direction is trivial. Suppose $\gf\ncong\bullet$. Since $\gf$ is rooted, we see that $\R[x]\neq\ve$ for any $x\in X$. Thus $\gf\md\D\top\vee\bd\top$. 
\end{proof}

\begin{proposition}\label{prop:axiom-L*-Kt}
    $L^*=\mathsf{K}_t\oplus(\D\top\vee\bd\top)$.
\end{proposition}
\begin{proof}
    Since $\bullet\not\md\D\top\vee\bd\top$ and $L^*=\mathsf{K}_t/\mathsf{Log}(\bullet)$, we get $L^*\sub\mathsf{K}\oplus(\D\top\vee\bd\top)$. Take any $\gf\in\mathsf{GF}_r(L^*)$. Then $\gf\ncong\bullet$. By Proposition~\ref{prop:deadend-frame}, $\gf\md\D\top\vee\bd\top$. Hence $\D\top\vee\bd\top\in\mathsf{Log}(\mathsf{GF}_r(L^*))=L^*$.
\end{proof}

\begin{proposition}\label{prop:splitting-Kt}
    Let $L\in\mathsf{NExt}(\mathsf{K}_t)$. Then $L\subsetneq L^*$ if and only if $L=\mathsf{K}_t$. 
\end{proposition}
\begin{proof}
    Suppose $L\subsetneq L^*$. Then $L\nsupseteq L^*$ and so $L\sub\mathsf{Log}(\bullet)$. Thus $\bullet\in\mathsf{Fr}(L)$. By Proposition~\ref{prop:deadend-frame}, $\mathsf{GF}_r(L)=\mathsf{GF}_r(\mathsf{K}_t)$, which entails $L\sub\mathsf{Log}(\mathsf{GF}_r(L))\sub\mathsf{Log}(\mathsf{GF}_r(\mathsf{K}_t))=\mathsf{K}_t$.
\end{proof}

\begin{lemma}\label{lem:no-more-iterated-splitting-Kt}
    There exists no logic splits $\NExt(L^*)$.
\end{lemma}
\begin{proof}
    Towards a contradiction, suppose there exists a logic $L$ splits $\NExt(L^*)$. Since $L\supseteq L^*$, $L\nsubseteq\mathsf{Log}(\bullet)$. By Theorem~\ref{thm:splittings-Kt-K4t}, $L$ does not split $\K$. By Fact~\ref{fact:splitting}, $\bigcap_{i\in I}L_i\sub L$ for a family of logics $\mathcal{K}=\cset{L_i\in\K:L_i\nsubseteq L, i\in I}$. Let $\mathcal{K}'=\cset{L_i\oplus(\D\top\vee\bd\top):i\in I}$. Take any $\phi\not\in L$. Since $L^*\sub L$, $(\D\top\vee\bd\top)\to\phi\not\in L$. Since $\bigcap_{i\in I}L_i\sub L$, $(\D\top\vee\bd\top)\to\phi\not\in L_i$ for some $i\in I$. Then there exists $\gf\in\mathsf{GF}(L_i)$ such that $\gf\not\md(\D\top\vee\bd\top)\to\phi$. By Proposition~\ref{prop:deadend-frame}, it is not hard to see that $\gf\md(\D\top\vee\bd\top)$ and $\gf\not\md\phi$. Thus $\phi\not\in L_i\oplus(\D\top\vee\bd\top)$. Since $\phi$ is arbitrarily chosen, we see that $\bigcap\mathcal{K}'\sub L$. Note that $L$ splits $\NExt(L^*)$ and $L\sub L_i\oplus L\nsubseteq L_1$ for any $i\in I$, we see that $\bigcap_{i\in I}L_i\oplus L\nsubseteq L_1$, which leads to a contradiction.
\end{proof}

\begin{theorem}\label{thm:iterated=union-Kt=L*}
    Let $L\in\K$. Then the following are equivalent:
    \begin{enumerate}[(1)]
        \item $L$ is an iterated splittings in $\K$.
        \item $L$ is a union-splitting in $\K$.
        \item $L\in\cset{\mathsf{K}_t,L^*}$.
    \end{enumerate}
\end{theorem}
\begin{proof}
    By Lemma~\ref{lem:no-more-iterated-splitting-Kt}, $\mathsf{K}_t$ and $L^*$ are the only two iterated splittings in $\K$. Thus (1) is equal to (3). By \cite[Theorem~22]{Kracht1992}, (2) is equal to (3).
\end{proof}

\begin{theorem}\label{thm:degreeL*-Kt}
    Let $L$ be a union-splitting in $\K$. Then $\mathsf{df}(L)=1$.
\end{theorem}
\begin{proof}
    By Theorem~\ref{thm:iterated=union-Kt=L*}, $L\in\cset{\mathsf{K}_t,L^*}$. Clearly, $\mathsf{df}_{\mathsf{K}_t}(\mathsf{K}_t)=1$. Let $L=L^*$. Suppose $\mathsf{Fin}(L')=\mathsf{Fin}(L^*)$ for some $L'\neq L^*$. Note that $L^*=\mathsf{K}_t\oplus(\D\top\vee\bd\top)$ has the FMP, we see $L'\subsetneq L^*$. By Proposition~\ref{prop:splitting-Kt}, $L'=\mathsf{K}_t$. Thus $\bullet\in\mathsf{Fin}(L')$, which contradicts to $\mathsf{Fin}(L')=\mathsf{Fin}(L^*)$.
\end{proof}

It is now sufficient to fix a tense logic $L\in\K\setminus\cset{\mathsf{K}_t,L^*}$ and show that $\mathsf{deg}(L)=2^{\aleph_0}$. The main proof idea is as follows: to show that there exists a continual family of tense logics which share the same Kripke frame as $L$, we construct for each $I\sub\mathbb{Z}^+$ a general frame $\gf_I$ such that $\mathsf{Fr}(L\cap\mathsf{Log}(\gf_I))=\mathsf{Fr}(L)$. Let $L_I=L\cap\mathsf{Log}(\gf_I)$. Then $\mathsf{deg}(L)\geq|\cset{L_I:I\sub\mathbb{Z}^+}|$ and we are done once we show that $L_I\neq L_J$ for any different $I,J\sub\mathbb{Z}^+$. 

Let $I\sub\mathbb{Z}^+$ be arbitrarily chosen and we start with the construction of $\gf_I$. Intuitively, the general frame $\gf_I$ will be a combination of a finite rooted frame $\F_L$ and a general frame $\gf'_I$. On one hand, the finite frame $\F_L$ is designed to refute some formula in $L$, since we have to ensure that $L_I\neq L_J$ for any different $I,J\sub\mathbb{Z}^+$, which requires $L\neq L\cap\mathsf{Log}(\gf_I)$. On the other hand, we have to construct $\gf'_I$ properly to make the logics $L_I$ pairwise different and share the same frames as $L$. The trick here is to choose $\F_L$ to be a finite frame of large enough r-degree, which ensures that $\gf'_I$ and $\F_L$ both work well after being combined. More precisely, we have the following lemma holds: 

\begin{lemma}\label{lem:large-enough-finite-model-Kt}
    For all $\phi\not\in L^*$ and $n\in\omega$, there exists $\F\in\mathsf{Fin}_r$ such that $\F\not\md\phi$ and $\mathsf{rdg}(\F)\geq n$.
\end{lemma}
\begin{proof}
    Take any $\phi\not\in L^*$. By Proposition~\ref{prop:axiom-L*-Kt}, $L^*=\mathsf{K}\oplus(\D\top\vee\bd\top)$ has the FMP and so $\G\not\md\phi$ for some $\G\in\mathsf{Fin}_r(L^*)$. By Propositions~\ref{prop:deadend-frame}, $\gf\ncong\bullet$. Then existence of required $\F$ follows from Corollary~\ref{coro:Kt-large-rdg}.
\end{proof}
By Proposition~\ref{prop:splitting-Kt}, there exists a formula $\phi_L\in L\setminus L^*$. By Lemma~\ref{lem:large-enough-finite-model-Kt}, there is a finite rooted frame $\F_L=(X_L,R_L)$ and $w_L,u_L\in X$ such that $\F_L,w_L\not\md\phi_L$ and $u_L\not\in R_\sharp^{\mathsf{md}(\phi)}[w_L]$.

\begin{definition}
    For each $I\sub\mathbb{Z}^+$, let $\gf'_I=(Y_I,S_I,B_I)$ be the general frame defined as follows:
\begin{itemize}
    \item $Y_I= \omega\cup\cset{i^*:i\in I\cup\cset{0}}$.
    \item $S_I= \cset{\tup{n,m}\in\omega\times\omega:n<m}\cup\cset{\tup{i^*,j}:i\in I\cup\cset{0}\text{ and }i\leq j}$.
    \item $B_I$ is the internal set generated by $\ve$.
\end{itemize}
\end{definition}
The general frames $\gf'_I$ were introduced in \cite{Chen.Ma2024}. The tense logic of $\gf'_I$ has no consistent proper extension and no Kripke frame \cite[Proposition 5.7]{Chen.Ma2024}.  
Now we define the general frame $\gf_I$ to be $(\tup{\F_Lu_L+0^*\kappa\gf'_I},A_I)$, where $A_I$ the internal set generated by $\mathcal{P}(X_L)$.

\begin{example}
    Let $\mathbb{P}$ be the set of all prime numbers. Then $\gf_\mathbb{P}$ is depicted by Figure \ref{fig:fi}.
    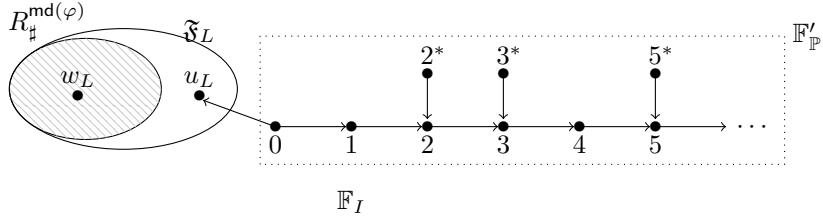
\begin{figure}[htbp]
    \[
    \begin{tikzpicture}
    \draw (-2,-1) node{$\gf_I$};
    \draw (-3,0) node{$\bullet$};
    \draw (-2,0) node{$\bullet$};
    \draw (0,0) node{$\bullet$};
    \draw (-1,0) node{$\bullet$};
    \draw (1,0) node{$\bullet$};
    \draw (2,0) node{$\bullet$};
    \draw (0,.7) node{$\bullet$};
    \draw (-1,.7) node{$\bullet$};
    \draw (2,.7) node{$\bullet$};
    \draw (-3,0) node[below]{$0$};
    \draw (-2,0) node[below]{$1$};
    \draw (-1,0) node[below]{$2$};
    \draw (0,0) node[below]{$3$};
    \draw (1,0) node[below]{$4$};
    \draw (2,0) node[below]{$5$};
    \draw (2.08,.7) node[above]{$5^\ast$};
    \draw (0.08,.7) node[above]{$3^\ast$};
    \draw (-.92,.7) node[above]{$2^\ast$};
    \draw (3.3,0) node{$\cdots$};
    \draw [->] (-3,0) -- (-2.07,0);
    \draw [->] (-2,0) -- (-1.07,0);
    \draw [->] (-1,0) -- (-0.07,0);
    \draw [->] (0,0) -- (0.93,0);
    \draw [->] (1,0) -- (1.93,0);
    \draw [->] (2,0) -- (2.93,0);
    \draw [<-] (0,0.1) -- (0,.655);
    \draw [<-] (-1,0.1) -- (-1,.655);
    \draw [<-] (2,0.1) -- (2,.655);
    
    \draw [dotted] (-3.2,-.5) rectangle (3.7,1.2);
    \draw (4,1.2) node{$\gf'_\mathbb{P}$};

    \draw (-4,1) node[above]{$\F_L$};
    \draw (-5.5,.5) ellipse (1 and .67);
    \fill [pattern =north west lines,pattern color = black!20](-5.5,.5) ellipse (1 and .67);
    \draw (-5,.5) ellipse (1.5 and .8);
    \draw (-4-2,1+0) node[above]{$R_\sharp^{\mathsf{md}(\phi)}$};
    \draw (-3.6-2,1+-.6) node{$\bullet$};
    \draw (-3.6-2,1+-.6) node[above]{$w_L$};
    \draw (-4,1+-.6) node{$\bullet$};
    \draw (-4,1+-.6) node[above]{$u_L$};
    \draw [<-] (-3.95,.35) -- (-3,0);
    \end{tikzpicture}
    \]
    \caption{The frame $\gf_\mathbb{P}$}
    \label{fig:fi}
    \end{figure}
\end{example}

Let $k\in\omega$ be such that $|\F_L|<k$ and $X_I={R_I}_\sharp^{k}[v]$ for all $v\in X_I$. For each $n\in\omega$ and $m\in\mathbb{Z}^+$, we define the formulas $\gamma_n$ and $\gamma^*_m$ as follows:
\begin{itemize}
    \item $\gamma_0=\bb\bot\wedge\D\bb^2\bot\wedge\D^k\bb^{k+1}\bot$ and $\gamma_{l+1}=\bd\gamma_l\wedge\bb^2\neg\gamma_l$.
    \item $\gamma^*_m=\D\gamma_m\wedge\B\neg\gamma_{m-1}\wedge\bb\bot\wedge\B\D^k\top$.
\end{itemize}

\begin{lemma}\label{lem:gamma-truth}
    For all $n\in\omega$, $m\in\mathbb{Z}^+$ and $x\in X_I$,
\begin{enumerate}[(1)]
    \item $\gf_I,x\md\gamma_n$ if and only if $x=n$;
    \item $\gf_I,x\md\gamma^*_m$ if and only if $m\in I$ and $x=m^*$.
\end{enumerate}
\end{lemma}
\begin{proof}
    For (1), we prove by induction on $n$. Let $n=0$. Note that $\breve{R_I}[0]=\ve$, $k\in R_I^k[0]$ and $\breve{R_I}^{k+1}[k]=\ve$, we have $\gf_I,0\md\gamma_0$. Suppose $\gf_I,x\md\gamma_0$. Then $\gf_I,x\md\D^k\bb^{k+1}\bot$, which entails that there exists a strict chain $\tup{x_i:i\leq k}$ with $x=x_0$. Since $|\F_L|<k$, we see $x\not\in X_L$. Since $\gf_I,x\md\bb\bot\wedge\D\bb^2\bot$, we have $x\not\in\mathbb{Z}^+\cup\cset{i^*:i\in I}$. Thus $x=0$. Let $n>0$. By induction hypothesis, for all $y\in X_I$, $\gf_I,y\md\gamma_{n-1}$ if and only if $y=n-1$. Since $n-1\in\breve{R_I}[n]\setminus\breve{R_I}^2[n]$, we have $\gf_I,n\md\gamma_n$. Suppose $\gf_I,x\md\gamma_n$. By induction hypothesis, $x\in R_I[n-1]\setminus R_I^2[n-1]$, which entails $x=n$.

    (2) The right-to-left direction is trivial. Suppose $\gf_I,x\md\gamma^*_m$. By $\gf_I,x\md\D\gamma_m\wedge\B\neg\gamma_{m-1}$, we have $x\not\in X_L\cup\cset{l,l^*\in X_I:l<m-1\text{ or }l>m}\cup\cset{m,(m-1)^*}$. Since $\gf_I,x\md\bb\bot\wedge\B\D^k\top$, $x\neq m-1$. By $x\in X_I=X_L\uplus(\omega\cup\cset{i^*:i\in I})$, we see $m\in I$ and $x=m^*$.
\end{proof}

\begin{lemma}\label{lem:no-frame-for-gfi-Kt}
    $\mathsf{Fr}(\mathsf{Log}(\gf_I))=\ve$.
\end{lemma}
\begin{proof}
    Suppose there exists $\G=(Y,S)\in\mathsf{Fr}(\mathsf{Log}(\gf_I))$. By an easy induction, we see that for all $V\in A_I$, $V\cap\omega$ and $V\cap\cset{i^*:i\in I\cup\cset{0}}$ are either finite or co-finite. Thus $\gf_I,1\md\B(\B(p\to\B p)\to p)\to\B p$. By Lemma~\ref{lem:gamma-truth}(1), $\gf_I\md\gamma_i\to\D\gamma_{i+1}$ for all $i\in\omega$. Thus we have
    \begin{center}
        $\G\md\gamma_1\to(\B(\B(p\to\B p)\to p)\to\B p)$ and $\G\md\cset{\gamma_i\to\D\gamma_{i+1}:i\in\omega}$.
    \end{center}
    Note that $\gf_I\md\Delta^k\gamma_1$ and $\gamma_1$ is variable-free, there exists $y\in Y$ such that $\G,y\md\gamma_1$. Since $\G\md\cset{\gamma_i\to\D\gamma_{i+1}:i\in\omega}$, we see that there exists an infinite strictly ascending $S$-chain $\tup{u_i:i\in\mathbb{Z}^+}$ such that $y=u_1$ and $\G,u_i\md\gamma_i$ for all $i\in\mathbb{Z}^+$. Let $U$ be a valuation in $\G$ such that $U(p)=\cset{u_{2i}:i\in\mathbb{Z}^+}$. Then we see that $\G,U,y\not\md\B(\B(p\to\B p)\to p)\to\B p$. Hence $\G\not\md\gamma_1\to(\B(\B(p\to\B p)\to p)\to\B p)$, which contradicts $\G\in\mathsf{Fr}(\mathsf{Log}(\gf_I))$.
\end{proof}

\begin{lemma}\label{lem:same-frame-Kt}
    For all $I\in\mathcal{P}(\mathbb{Z}^+)$, $\mathsf{Fr}(L)=\mathsf{Fr}(L_I)$.
\end{lemma}
\begin{proof}
    By Lemmas~\ref{lem:intersection-rooted-frame} and \ref{lem:no-frame-for-gfi-Kt}, $\mathsf{Fr}_r(L_I)=\mathsf{Fr}_r(L)\cup\mathsf{Fr}_r(\mathsf{Log}(\gf_I))=\mathsf{Fr}_r(L)$.
\end{proof}

\begin{lemma}\label{lem:pairwise-diff-Kt}
    For all $I,J\in\mathcal{P}(\mathbb{Z}^+)$, $I\neq J$ implies $L_I\neq L_J$.
\end{lemma}
\begin{proof}
    Take any distinct $I,J\in\mathcal{P}(\mathbb{Z}^+)$. Let $i\in I\setminus J$. It suffices to show that 
    \begin{center}
        $\neg\phi_L\to\Delta^{k}\gamma^*_{i}\in L_I\setminus L_J$.
    \end{center}
    By Lemma \ref{lem:gamma-truth}(2), $\gf_I,i^*\md\gamma^*_i$ and so $\gf_I\md\Delta^k\gamma^*_i$. Since $\phi_L\in L$, we see that $\neg\phi_L\to\Delta^{k}\gamma^*_{i}\in L\cap\mathsf{Log}(\gf_I)= L_I$. Since $i\not\in J$, by Lemma \ref{lem:gamma-truth}(2), we see $\gf_J\md\neg\gamma^*_i$. Note that $\gf_J\rsto{R_J}_\sharp^{\mathsf{md}(\phi_L)}[w_L]\cong\F_L\rsto{R_L}_\sharp^{\mathsf{md}(\phi_L)}[w_L]$, we have $\gf_J,w_L\not\md\phi_L$. Thus $\gf_J,w_L\not\md\phi_L\vee\Delta^k\gamma^*_i$ and so $\neg\phi_L\to\Delta^{k}\gamma^*_{i}\not\in L_J$.
\end{proof}

Since $L$ is chosen arbitrarily, by Lemmas~\ref{lem:same-frame-Kt} and \ref{lem:pairwise-diff-Kt}, we have
\begin{theorem}\label{thm:degreeKt-continuum}
    Let $L\in\K\setminus\cset{\mathsf{K}_t,L^*}$. Then $\mathsf{deg}(L)=2^{\aleph_0}$.
\end{theorem}

We conclude this section by the following characterization theorems in $\mathsf{NExt}(\mathsf{K}_t)$:
\begin{theorem}\label{thm:degreeKt-main}
    Let $L\in\K$. Then the following are equivalent:
    \begin{enumerate}[(1)]
        \item $L$ is a union-splitting in $\K$.
        \item $L$ is an iterated splitting in $\K$.
        \item $\mathsf{df}(L)=1$.
        \item $\mathsf{deg}(L)=1$.
        \item $\mathsf{df}(L)\neq 2^{\aleph_0}$.
        \item $\mathsf{deg}(L)\neq 2^{\aleph_0}$.
    \end{enumerate}    
\end{theorem}
\begin{proof}
    The equivalence of (1), (2) and (3) follows from Theorem~\ref{thm:iterated=union-Kt=L*}. Since $|\K|\leq 2^{\aleph_0}$, by Proposition~\ref{prop:deg-leq-df}, (3) implies (4), and (5) implies (6). Clearly (3) implies (5), and (4) implies (6). Finally, by Theorem~\ref{thm:degreeKt-continuum}, (6) implies (1), which concludes the proof.
\end{proof}

The following dichotomy theorem follows immediately:
\begin{theorem}\label{thm:dichotomy-Kt}
    For all $L\in\K$, $\mathsf{deg}(L)=\mathsf{df}(L)\in\cset{1,2^{\aleph_0}}$.
\end{theorem}

\subsection{Degree of Kripke-incompleteness in $\LT$}\label{subsec:degree-K4t}
The proof idea of the dichotomy theorem for $\LT$ is similar to the one for $\K$. In this section, we always presume that frames are transitive and write $\mathsf{deg}$ for $\mathsf{deg}_{\mathsf{K4}_t}$. By \cite[Theorem 10]{Kracht1992}, $\tup{\mathsf{Log}(\bullet),\mathsf{K4}_t/\bullet}$ is the unique splitting pair in $\LT$. In this section, we write $L^*$ for $\mathsf{K4}_t/\bullet$. Similar to the case for $\K$, the logics $\mathsf{K4}_t$ and $L^*$ are exactly the union-splittings in $\LT$. Since $\mathsf{K4}_t$ is Kripke complete, $\mathsf{deg}(\mathsf{K4}_t)=1$. Moreover, by arguments similar to those in Section~\ref{subsec:degree-Kt}, we have

\begin{proposition}\label{prop:L*-K4t}
    Let $\gf\in\mathsf{GF}_r(\mathsf{K4}_t)$. Then (i) $\gf\md\D\top\vee\bd\top$ if and only if $\gf\ncong\bullet$; (ii) $L^*=\mathsf{K4}\oplus(\D\top\vee\bd\top)$; and (iii) for all $L\in\mathsf{NExt}(\mathsf{K4}_t)$, $L\subsetneq L^*$ if and only if $L=\mathsf{K4}_t$.
\end{proposition}

\begin{theorem}\label{thm:iterated=union-K4t=L*}
    Let $L\in\LT$. Then the following are equivalent:
    \begin{enumerate}[(1)]
        \item $L$ is an iterated splittings in $\LT$.
        \item $L$ is a union-splitting in $\LT$.
        \item $L\in\cset{\mathsf{K4}_t,L^*}$.
    \end{enumerate}
\end{theorem}

\begin{theorem}\label{thm:degreeK4tL*}
    For all union-splittings $L$ in $\LT$, $\mathsf{df}_{\mathsf{K4}_t}(L)=\mathsf{deg}_{\mathsf{K4}_t}(L)=1$.
\end{theorem}

It is sufficient now to prove that $\mathsf{deg}(L)=2^{\aleph_0}$ for all $L\in\LT\setminus\cset{\mathsf{K4}_t,L^*}$. The method we use here is similar to the one in Section~\ref{subsec:degree-Kt} and the key lemma is the following:

\begin{lemma}\label{lem:large-enough-finite-model}
    Let $\phi\not\in L^*$. Then (i) $\phi$ is refuted by some non-symmetric $\G\in\mathsf{Fin}_r(\mathsf{K4}_t)$; (ii) for all $n\in\omega$, there exists $\F\in\mathsf{Fin}_r(\mathsf{K4}_t)$ such that $\F\not\md\phi$ and $\mathsf{rdg}(\F)\geq n$.
\end{lemma}
\begin{proof}
    For (i), take any $\phi\not\in L^*$. By Proposition~\ref{prop:L*-K4t}, $L^*=\mathsf{K4}_t\oplus(\D\top\vee\bd\top)$ has the FMP. Then there exists $\F'=(X',R')\in\mathsf{Fin}_r(L^*)$ such that $\F'\not\md\phi$. If $\F'$ is already non-symmetric, then take $\G=\F'$ we are done. Suppose $\F'$ is symmetric. Since $\F'$ is rooted and transitive, $\F$ is a cluster. Note that $\F'\md\D\top\vee\bd\top$, we see that $R'_\sharp[x']\neq\ve$ for each $x'\in X'$ and so $R'=X'\times X'$. Let $\G=(X,R)$ where $X=X'\times\cset{0,1}$ and $R=\cset{\tup{\tup{x,a},\tup{y,b}}\in X\times X:a\leq b}$. Consider the map $f:X\to X'$ defined by $f:\tup{x,a}\mapsto x$ for all $\tup{x,a}\in X$. Obviously, $f$ is a t-morphism from $\G$ to $\F'$. Thus $\G\not\md\phi$. Hence (i) holds. By Corollary~\ref{coro:Tr-large-rdg}, (ii) follows from (i) immediately.
\end{proof}

We can now prove the following theorem:

\begin{theorem}\label{thm:degreeK4t-continuum}
    Let $L\in\LT\setminus\cset{\mathsf{K4}_t,L^*}$. Then $\mathsf{deg}(L)=2^{\aleph_0}$.
\end{theorem}
\begin{proof}
    Let $L\in\LT\setminus\cset{\mathsf{K4}_t,L^*}$. By Proposition~\ref{prop:L*-K4t}, $L\nsubseteq L^*$ and so there exists $\phi_L\in L\setminus L^*$. By Lemma \ref{lem:large-enough-finite-model}, there is a finite rooted frame $\F_L=(X_L,R_L)$ and $w_L,u_L\in X$ such that $\F_L,w_L\not\md\phi_L$ and $u_L\not\in R_\sharp^{\mathsf{md}(\phi)}[w_L]$. For each $I\in\mathbb{Z}^+$, we define $\gf_I$ to be the general frame $(\tup{\F_Lu_L+^t0^*\kappa\gf'_I},A_I)$, where $A_I$ the internal set generated by $\mathcal{P}(X_L)$. Then we see that $\gf$ is transitive. By almost same proofs as those of Lemmas~\ref{lem:gamma-truth}, \ref{lem:no-frame-for-gfi-Kt}, \ref{lem:same-frame-Kt} and \ref{lem:pairwise-diff-Kt}, we have
    \begin{itemize}
        \item $\mathsf{Fr}(L)=\mathsf{Fr}(L_I)$ for all $I\sub\mathbb{Z}^+$;
        \item $I\neq J$ implies $L_I\neq L_J$ for all $I,J\sub\mathbb{Z}^+$.
    \end{itemize}
    It follows immediately that $2^{\aleph_0}=|\cset{L_I:I\sub\mathbb{Z}^+}|\leq\mathsf{deg}(L)\leq 2^{\aleph_0}$. Hence $\mathsf{deg}(L)=2^{\aleph_0}$.
\end{proof}

By Theorems~\ref{thm:iterated=union-K4t=L*}, \ref{thm:degreeK4tL*} and \ref{thm:degreeK4t-continuum}, we have the following characterization theorem:
\begin{theorem}\label{thm:degreeK4t-main}
    Let $L\in\LT$. Then the following are equivalent:
    \begin{enumerate}[(1)]
        \item $L$ is a union-splitting in $\LT$.
        \item $L$ is an iterated splitting in $\LT$.
        \item $\mathsf{df}(L)=1$.
        \item $\mathsf{deg}(L)=1$.
        \item $\mathsf{df}(L)\neq 2^{\aleph_0}$.
        \item $\mathsf{deg}(L)\neq 2^{\aleph_0}$.
    \end{enumerate}    
\end{theorem}

Again, the following dichotomy theorem follows immediately:
\begin{theorem}\label{thm:dichotomy-K4t}
    For all $L\in\LT$, $\mathsf{deg}(L)=\mathsf{df}(L)\in\cset{1,2^{\aleph_0}}$.
\end{theorem}

\section{Kripke-incompleteness in $\mathsf{NExt}(\mathsf{S4}_t)$}\label{sec:DKI-S4t}

In this section, we focus on the degree of Kripke-incompleteness in $\NExt(\ST)$. In what follows, we write $\mathsf{deg}(L)$ for the degree of Kripke-incompleteness of $L$ in $\NExt(\ST)$. Our aim is to prove the dichotomy theorem for $\NExt(\ST)$. We first characterize the iterated splittings in $\NExt(\ST)$ and show that every iterated splitting has the degree of Kripke-incompleteness $1$. Then by constructing some frames generalizing the Nishimura-Rieger ladder, we show that all other tense logics in $\NExt(\ST)$ are of the degree of Kripke-incompleteness $2^{\aleph_0}$, which gives a characterization theorem for the degree of Kripke-incompleteness in $\NExt(\ST)$ and completes the proof of the dichotomy theorem. 

Let $n\in\mathbb{Z}^+$. Then we write $\C_n$ for the reflexive transitive chain of length $n$ and $\mathfrak{Cl}$ for the $n$-cluster, that is, $\C_n=(n,\leq)$ and $\mathfrak{Cl}_n=(n,n\times n)$.

\begin{lemma}\label{lem:iterated-splitting-S4t}
    Let $L\in\NExt(\ST)$. Then $L$ is an iterated splitting if and only if $L\in\NExt(\mathsf{S5}_t)\cup\cset{\ST}$.
\end{lemma}
\begin{proof}
    For the right-to-left direction, note first that $\ST$ and $\mathsf{S5}_t=\ST/\C_2$ are iterated splittings. By \cite[Lemma 6.6]{Chen2024}, $\NExt(\mathsf{S5}_t)$ is isomorphic to $(\omega,\geq)$ and every consistent extension of $\mathsf{S5}_t$ is of the form $\mathsf{Log}(\mathfrak{Cl}_k)$. Then we have $\L=\mathsf{S5}_t/\mathsf{Log}(\mathfrak{Cl}_1)$ and $\mathsf{Log}(\mathfrak{Cl}_k)=\mathsf{S5}_t/\mathsf{Log}(\mathfrak{Cl}_{k+1})$ for all $k\in\mathbb{Z}^+$. Suppose $L=\ST/L_1/\cdots/L_n$ is an iterated splitting. If $n=0$, then $L=\ST$. If $n=1$, by \cite[Theorem 21]{Kracht1992}, $L\in\cset{\mathsf{S5}_t,\L}\sub\NExt(\mathsf{S5}_t)$. If $n\geq 2$, then $L\supseteq L/L_1\supseteq\mathsf{S5}_t$. Thus $L\in\NExt(\mathsf{S5}_t)\cup\cset{\ST}$.
\end{proof}

\begin{lemma}\label{lem:deg=1}
    Let $L$ be an iterated splitting in $\NExt(\ST)$. Then $\mathsf{df}(L)=1$.
\end{lemma}
\begin{proof}
    By Lemma~\ref{lem:iterated-splitting-S4t}, $L\in\NExt(\mathsf{S5}_t)\cup\cset{\ST}$. Clearly, $\mathsf{df}(\ST)=1$. Let $L\in\NExt(\mathsf{S5}_t)$. Take any $L'\in\NExt(\ST)$ such that $\mathsf{Fin}(L')=\mathsf{Fin}(L)$. Then $\C_2\not\md L'$. By \cite[Theorem 21]{Kracht1992}, $\tup{\mathsf{Log}(\C_2),\mathsf{S5}_t}$ is a splitting pair in $\NExt(\ST)$, which entails $L'\in\NExt(\mathsf{S5}_t)$. By \cite[Lemma 6.6]{Chen2024}, every extension of $\mathsf{S5}_t$ has the FMP and so $L=\mathsf{Log}(\mathsf{Fin}(L))=\mathsf{Log}(\mathsf{Fin}(L'))=L'$. Hence $\mathsf{df}(L)=1$.
\end{proof}
Recall from \cite[Theorem 21]{Kracht1992} that $\tup{\mathsf{Log}(\C_2),\mathsf{S5}_t}$ and $\tup{\mathsf{Log}(\C_1),\L}$ are the only two splitting pairs in $\NExt(\ST)$. Then $\cset{\L,\mathsf{S5}_t,\ST}$ is exactly the set of union-splittings in $\NExt(\ST)$ and so every union-splitting has the degree of Kripke-incompleteness $1$. However, since there are countably many tense logics in $\NExt(\mathsf{S5}_t)$, Lemma~\ref{lem:deg=1} indicates that logics of the degree of Kripke-incompleteness $1$ are not necessary union-splittings, which shows that Blok's characterization of the degree of Kripke-incompleteness for $\NExt(\mathsf{K})$ can not be generalized to $\NExt(\ST)$.

\begin{lemma}\label{lem:large-enough-finite-model-S4t}
    Let $\phi\not\in\ST$. Then (i) $\phi$ is refuted by some non-symmetric $\G\in\mathsf{Fin}_r(\mathsf{S4}_t)$; (ii) for all $n\in\omega$, there exists $\F\in\mathsf{Fin}_r(\mathsf{S4}_t)$ such that $\F\not\md\phi$ and $\mathsf{rdg}(\F)\geq n$.
\end{lemma}
\begin{proof}
    For (i), take any $\phi\not\in\ST$. Since $\ST$ has the FMP, there exists $\F'=(X',R')\in\mathsf{Fin}_r(\ST)$ such that $\F'\not\md\phi$. If $\F'$ is non-symmetric, then take $\G=\F'$ we are done. Suppose $\F'$ is symmetric. Then $R'=X'\times X'$. Let $\G=(X,R)$ where $X=X'\times\cset{0,1}$ and $R=\cset{\tup{\tup{x,a},\tup{y,b}}\in X\times X:a\leq b}$. Then the map $f:X\to X'$ defined by $f:\tup{x,a}\mapsto x$ for all $\tup{x,a}\in X$ is a t-morphism from $\G$ to $\F'$. Thus $\G\not\md\phi$. Hence (i) holds. By Corollary~\ref{coro:Tr-large-rdg}, (ii) follows from (i) immediately.
\end{proof}
Let $L\in\NExt(\ST)$ be an arbitrarily fixed logic such that $L\not\in\NExt(\mathsf{S5}_t)\cup\cset{\ST}$. Since $L\supsetneq\ST$, there exists $\phi_L\in L\setminus\mathsf{S4}_t$. Then there exists $\F_L\in\mathsf{Fin}_r(\ST)$ and $w_L,u_L\in X_L$ with $\F_L,w_L\not\md\phi_L$ and $u_L\not\in R_\sharp^{\mathsf{md}(\phi)}[w_L]$.
As what we did in Section~\ref{sec:DKI-Kt-K4t}, to show that $L$ is of the degree of Kripke-incompleteness $2^{\aleph_0}$, it suffices to construct a continual family of general frames $\tup{\gf_I:I\sub\mathbb{Z}^+}$ such that $\Fr(L)=\Fr(L\cap\mathsf{Log}(\gf_I))$ and $L\cap\mathsf{Log}(\gf_I)\neq L\cap\mathsf{Log}(\gf_I)$ for any $I\neq J\sub\mathbb{Z}^+$.

\begin{definition}
    Let $I\in\mathcal{P}(\mathbb{Z}^+)$. The frame $\F'_I=(Y_I,S_I)$ is defined as follows:
    \begin{itemize}
        \item $Y_I=A\cup B\cup C_I\cup\cset{x_0,x_1,x_2,y_0,y_1,r_0,r_1,r_2,r'}$, where $A=\cset{a_i:i\in\omega}$, $B=\cset{b_i:i\in\omega}$ and $C=\cset{c_i:i\in I\cup\cset{0}}$;
        \item $S_I$ is the reflexive-transitive closure of the union of the following binary relations:
        \begin{itemize}
            \item $\cset{\tup{x_0,x_1},\tup{x_2,x_1},\tup{x_2,a_0},\tup{y_1,y_0},\tup{y_1,b_0},\tup{r_0,r'},\tup{r_0,r_1},\tup{r_2,r_1}}$;
            \item $\cset{\tup{c_i,c_j}:i>j\text{ and }i,j\in I\cup\cset{0}}$;
            \item $\cset{\tup{a_i,a_j}:i>j\in\omega}\cup\cset{\tup{a_i,b_j}:i>j\in\omega}$;
            \item $\cset{\tup{b_i,b_j}:i>j\in\omega}\cup\cset{\tup{b_i,a_j}:i>j+1\in\omega}$.
        \end{itemize}
    \end{itemize}
    The general frame $\gf_I=(X_I,R_I,A_I)$, where $\F_I=(X_I,R_I)=\tup{\F_Lu_L+^tr_3\F'_I}$ and $A_I$ is the internal set generated by $\mathcal{P}(X_L)$.
\end{definition}
An example of the underlying frame $\F_I$ is as depicted in Figure~\ref{fig:fi-S4t}. Clearly, $\gf_I\in\mathsf{GF}(\ST)$. The reader might notice that Rieger-Nishimura ladder can be embedded into $\F'_I$ by an order-preserving map. Let us take a closer look at $\gf_I$ by showing some of its properties.

\begin{figure}[htbp]
    \small
\[
\begin{tikzpicture}[scale=0.9]
\draw (-4,5) node{$\gf_I$};

\def\ptRad{.2pt}
\node (a0) at (0,8)[label=left:$a_0$]{$\circ$};
\node (a1) at (0,7)[label=left:$a_1$]{$\circ$};
\node (a2) at (0,6)[label=left:$a_2$]{$\circ$};
\node (a3) at (0,5)[label=left:$a_3$]{$\circ$};
\node (a4) at (0,4)[label=left:$a_4$]{$\circ$};
\node (b0) at (1.5,8)[label=right:$b_0$]{$\circ$};
\node (b1) at (1.5,7)[label=right:$b_1$]{$\circ$};
\node (b2) at (1.50,6)[label=right:$b_2$]{$\circ$};
\node (b3) at (1.5,5)[label=right:$b_3$]{$\circ$};
\node (b4) at (1.5,4)[label=right:$b_4$]{$\circ$};   
\node (ai0) at (0,2.5)[label=left:$a_{i-1}$]{$\circ$};
\node (ai1) at (0,1.5)[label=left:$a_{i}$]{$\circ$};
\node (ai2) at (0,0.5)[label=left:$a_{i+1}$]{$\circ$};
\node (bi0) at (1.5,2.5)[label=right:$b_{i-1}$]{$\circ$};
\node (bi1) at (1.5,1.5)[label=right:$b_{i}$]{$\circ$};
\node (bi2) at (1.50,0.5)[label=right:$b_{i+1}$]{$\circ$};
\node (r) at (.75,-1)[label=below:$r_0$]{$\circ$};
\node (r3) at (2.5,-1)[label=right:$r'$]{$\circ$};
\node (r1) at (-2,-1)[label=left:$r_1$]{$\circ$};
\node (vd1) at (.75,3.4){$\vdots$};
\node (vd2) at (0,3.4){$\vdots$};
\node (vd3) at (1.5,3.4){$\vdots$};
\node (vd4) at (0,0){$\vdots$};
\node (vd4) at (.75,0){$\vdots$};
\node (vd4) at (1.5,0){$\vdots$};
\node (vd4) at (-1.25,0){$\vdots$};

\node (x0) at (-3,7)[label=left:$x_0$]{$\circ$};
\node (x1) at (-2,8)[label=left:$x_1$]{$\circ$};
\node (x2) at (-1,7)[label=left:$x_2$]{$\circ$};
\node (y0) at (3.5,8)[label=right:$y_0$]{$\circ$};
\node (y1) at (2.5,7)[label=right:$y_1$]{$\circ$};

\node (z20) at (-1.75,6)[label=left:$c_0$]{$\circ$};
\node (z21) at (-1.25,5)[label=left:$c_2$]{$\circ$};
\node (z31) at (-1.25,4)[label=left:$c_3$,label=below:$\vdots$]{$\circ$};
\node (zi1) at (-1.25,.5)[label=left:$c_i$,label=above:$\vdots$]{$\circ$};

\draw [->,shorten <>=\ptRad] (a1) -- (a0);
\draw [->,shorten <>=\ptRad] (a2) -- (a1);
\draw [->,shorten <>=\ptRad] (a3) -- (a2);
\draw [->,shorten <>=\ptRad] (a4) -- (a3);
\draw [->,shorten <>=\ptRad] (b1) -- (b0);
\draw [->,shorten <>=\ptRad] (b2) -- (b1);
\draw [->,shorten <>=\ptRad] (b3) -- (b2);
\draw [->,shorten <>=\ptRad] (b4) -- (b3);
\draw [->,shorten <>=\ptRad] (a1) -- (b0);
\draw [->,shorten <>=\ptRad] (a2) -- (b1);
\draw [->,shorten <>=\ptRad] (a3) -- (b2);
\draw [->,shorten <>=\ptRad] (a4) -- (b3);
\draw [->,shorten <>=\ptRad] (b2) -- (a0);
\draw [->,shorten <>=\ptRad] (b3) -- (a1);
\draw [->,shorten <>=\ptRad] (b4) -- (a2);
\draw [->,shorten <>=\ptRad] (ai1) -- (ai0);
\draw [->,shorten <>=\ptRad] (ai2) -- (ai1);
\draw [->,shorten <>=\ptRad] (bi1) -- (bi0);
\draw [->,shorten <>=\ptRad] (bi2) -- (bi1);
\draw [->,shorten <>=\ptRad] (ai1) -- (bi0);
\draw [->,shorten <>=\ptRad] (ai2) -- (bi1);
\draw [->,shorten <>=\ptRad] (bi2) -- (ai0);
\draw [->,shorten <>=\ptRad] (x0) -- (x1);
\draw [->,shorten <>=\ptRad] (x2) -- (x1);
\draw [->,shorten <>=\ptRad] (x2) -- (a0);
\draw [->,shorten <>=\ptRad] (y1) -- (y0);
\draw [->,shorten <>=\ptRad] (y1) -- (b0);

\draw [->,shorten <>=\ptRad] (z21) -- (z20);
\draw [->,shorten <>=\ptRad] (z21) -- (a2);
\draw [->,shorten <>=\ptRad] (z31) -- (z21);
\draw [->,shorten <>=\ptRad] (z31) -- (a3);
\draw [->,shorten <>=\ptRad] (zi1) -- (ai1);

\draw [->,shorten <>=\ptRad] (r) -- (r1);
\draw [->,shorten <>=\ptRad] (r) -- (r3);
\draw [->,shorten <>=\ptRad] (r) -- (0,-.5);
\draw [->,shorten <>=\ptRad] (r) -- (1.5,-.5);
\draw [->,shorten <>=\ptRad] (r) -- (-1.25,-.5);

\node (ul) at (-4,.4)[label=above:$u_L$]{$\circ$};
\draw (-4,1.2) node[above]{$\F_L$};
\draw (-5.5,.5) ellipse (1 and .67);
\fill [pattern =north west lines,pattern color = black!20](-5.5,.5) ellipse (1 and .67);
\draw (-5,.5) ellipse (1.5 and .8);
\draw (-4-2,1.3+0) node[above]{$R_\sharp^{\mathsf{md}(\phi)}[w_L]$};
\draw (-3.6-2,1+-.6) node{$\circ$};
\draw (-3.6-2,1+-.6) node[above]{$w_L$};
\draw [->,shorten <>=\ptRad] (ul) -- (r1);
\end{tikzpicture}
\]
\caption{The frame $\F_I$ where $1\not\in I$ and $2,3,i\in I$}
\label{fig:fi-S4t}
\end{figure}

\begin{lemma}\label{lem:gfi-S4t-prop}
    Let $k\in\omega$ be such that $|\F_L|+5<k$ and $X_I={R_I}_\sharp^{k}[v]$ for all $v\in X_I$. Then
    \begin{enumerate}[(1)]
        \item $\gf_I\md\mathsf{bw}^+_k\wedge\mathsf{bw}^-_k\wedge\mathsf{bz}_k$;
        \item $\gf_I\md(\mathsf{grz}^+\wedge\mathsf{grz}^-)\vee(\mathsf{alt}^+_k\wedge\mathsf{alt}^-_k)$.
    \end{enumerate}
\end{lemma}
\begin{proof}
    For (1), by Fact~\ref{fact:bounds}, it suffices to notice that there exists no anti-chain or zigzag of size greater than $k$.
    For (2), by Fact~\ref{fact:bounds}, $\gf_I,x\md\mathsf{alt}^+_k\wedge\mathsf{alt}^-_k$ for all $x\in X_L\cup\cset{r_1}$ and $\gf_I,x\md\mathsf{grz}^+$ for all $x\not\in X_L\cup\cset{r_1}$. We now claim that the following holds:
    \begin{center}
        (\dag) For all admissible set $U\in A_I$ and chain $C\sub A\cup B\cup C_I$, either $U\cap C$ of $C\setminus U$ is finite.
    \end{center}
    We prove (\dag) by induction on the construction of $U$. Clearly, either $U\cap C$ of $C\setminus U$ is finite for all $U\in\mathcal{P}(X_I)$. The Boolean cases are straightforward. Let $U=\breve{R}[V]$. If $C\cap V\neq\ve$, then $C\setminus U$ is finite. If $C\cap V=\ve$, then $C\cap U=\ve$ is finite. Let $U=R[V]$. Take any chain $C\sub A\cup B\cup C_I$. If $r_0\in V\cap C$, then $C\setminus U=\ve$ is finite. Suppose $r_0\not\in V\cap C$. By induction hypothesis, $C\cap V$ or $C\setminus V$ is finite. Note that $C\cap U$ and $C\setminus U$ is finite iff $C\cap V$ and $C\setminus V$ is finite, respectively. Thus (\dag) holds.
    Take any $x_0\not\in X_L\cup\cset{r_1}$. Suppose $\gf_I,V,x_0\not\md\mathsf{grz}^-$ for some valuation $V$ in $\gf_I$. Then $x_0\not\in V(p)$ and $x_1\in V(\neg p\wedge\bd p)$ for some $x_1\in\breve{R}[x_0]$. By repeating this construction, there exists a co-chain $C=\cset{x_i:i\in\omega}\sub\breve{R}[x]$ such that $C\cap V(\neg p)=\cset{x_{2i}:i\in\omega}$. Thus $|C\cap V(p)|=|C\setminus V(p)|=\aleph_0$, which contradicts (\dag). Hence $\gf_I,x\md\mathsf{grz}^-$ and so $\gf_I\md(\mathsf{grz}^+\wedge\mathsf{grz}^-)\vee(\mathsf{alt}^+_k\wedge\mathsf{alt}^-_k)$.
\end{proof}

\begin{lemma}\label{lem:no-inf-Kripke-frame-S4tFi}
    $\mathsf{Fin}_r(\mathsf{Log}(\gf_I))=\mathsf{Fr}_r(\mathsf{Log}(\gf_I))$.
\end{lemma}
\begin{proof}
    Take any $\G=(Y,S)\in\mathsf{Fr}_r(\mathsf{Log}(\gf_I))$. By Lemma~\ref{lem:gfi-S4t-prop}(1), $\G$ is $k$-transitive and of both forth-width and back-width no more than $k$. By Lemma~\ref{lem:gfi-S4t-prop}(2) and Fact~\ref{fact:grz}, there is no infinite chain or infinite cluster in $\G$. Then the readers can readily check that $S_\sharp[y]$ is finite for all $y\in Y$, which entails $\G\in\mathsf{Fin}_r(\mathsf{Log}(\gf_I))$.
\end{proof}

Let $L_I=L\cap\mathsf{Log}(\gf_I)$. To show $\Fr(L)=\Fr(L_I)$, we introduce some auxiliary notions from \cite{Chen2024}.

\begin{lemma}\label{lem:t-morphism-source}
    Let $\gf=(X,R,A)\in\mathsf{GF}$, $\G=(Y,S)\in\Fr$, $f:\gf\twoheadrightarrow\G$ and $x,y\in X$. Suppose $f(x)=f(y)$. Then $f[R[x]]=f[R[y]]$ and $f[\breve{R}[x]]=f[\breve{R}[y]]$. %
\end{lemma}
\begin{proof}
    By $f:\gf\twoheadrightarrow\G$, $f[R[x]]=S[f(x)]=S[f(y)]=f[R[y]]$. Symmetrically, $f[\breve{R}[x]]=f[\breve{R}[y]]$. 
\end{proof}

\begin{definition}
    Let $\gf=(X,R,A)\in\mathsf{GF}_r$, $\G=(Y,S)\in\Fr_r$ and $f:\gf\twoheadrightarrow\G$. A subset $Z\sub X$ is called {\em sufficient} if for all $z\in Z$, there exist $u,v\in Z$ such that $f(z)=f(u)=f(v)$ and $R[u]\cup\breve{R}[v]\sub Z$. 
\end{definition}

\begin{lemma}\label{lem:sufficient-full}
    Let $\gf=(X,R,A)\in\mathsf{GF}_r$, $\G=(Y,S)\in\Fr_r$ and $f:\gf\twoheadrightarrow\G$. Suppose $Z\sub X$ is sufficient. Then $f[Z]=Y$.
\end{lemma}
\begin{proof}
    Take any $z\in Z$. It suffices to show $S_\sharp^n[f(z)]\sub f[Z]$ for all $n\in\omega$. The proof proceeds by induction on $n$. For details, see the proof of \cite[Lemma 3.8]{Chen2024}.
\end{proof}

It is no longer possible to show that $\mathsf{Log}(\gf_I)$ has no Kripke frame as in Section~\ref{sec:DKI-Kt-K4t}. However, we can still show that adding the general frame $\gf_I$ brings no new Kripke frame to $L$.

\begin{lemma}\label{lem:no-nontrivial-image}
    $\Fr_r(\mathsf{Log}(\gf_I))=\mathsf{TM}(\C_2)$.
\end{lemma}
\begin{proof}
    Note that $\mathsf{Log}(\gf_I)\not\in\NExt(\mathsf{S5}_t)$ and $\tup{\mathsf{Log}(\C_2),\mathsf{S5}_t}$ is a splitting pair in $\NExt(\ST)$, we have $\C_2\md\mathsf{Log}(\gf_I)$ and so $\Fr_r(\mathsf{Log}(\gf_I))\supseteq\mathsf{TM}(\C_2)$. Take any $\G\in\Fr_r(\mathsf{Log}(\gf_I))$. By Lemma~\ref{lem:no-inf-Kripke-frame-S4tFi}, $\G$ is finite. Let $\J^{k}(\G)$ be the Jankov-formula of $\G$ of degree $k$. Since $\G\not\md\neg\J^{k}(\G)$, we see that $\gf_I\not\md\neg\J^{k}(\G)$. Note that $\gf_I$ is $k$-transitive and rooted, by Theorem~\ref{thm:JankovLemma}, $\G$ is a t-morphic image of $\gf_I$. Let $f:\gf_I\twoheadrightarrow\G$. 
    Now it suffices to show that $\G\in\mathsf{TM}(\C_2)$. Suppose $\G\not\in\mathsf{TM}(\C_2)$. Then

    \noindent\textbf{Claim 1.} For all $u\in X_I\setminus X_L$, $f(u)=f(x_0)$ implies $u=x_0$.

    \noindent\textbf{Proof of Claim 1.} Take any $u\in X_I\setminus X_L$ such that $f(u)=f(x_0)$. Towards a contradiction, suppose also $u\neq x_0$. Then we have the following cases:
    \begin{enumerate}[(1)]
        \item $u\in\cset{x_1,y_0,r'}$. By Lemma~\ref{lem:t-morphism-source}, $S[f(x_0)]=S[f(u)]=f[R[u]]=\cset{f(u)}$ and $\breve{S}[f(x_0)]=f[\breve{R}[x_0]]=\cset{f(x_0)}$. Thus $\G\iso\C_1$ and contradicts the assumption.
        \item $u=x_2$. Then clearly $\cset{x_0,x_1,x_2}$ is sufficient, which entails $Y=\cset{f(x_0),f(x_1)}$. Note that $f(x_0)\not\in S[f(x_1)]$, $\G\cong\C_2$ and contradicts the assumption.
        \item $u=y_1$. Similar to (2), we see $\cset{x_0,x_1,y_0,y_1}$ is sufficient and so $\G\cong\C_2$, which contradicts the assumption.
        \item None of (1)-(3) holds. Then $u\in R[r_0]\setminus\cset{r'}$. By Lemma~\ref{lem:t-morphism-source}, $f(r_0)\in f[\breve{R}[u]]=f[\breve{R}[x_0]]=\cset{f(x_0)}$. By Lemma~\ref{lem:t-morphism-source} again, we see that $f(r')\in\cset{f(x_0),f(x_1)}$. By (1), $f(r')=f(x_1)$, which entails that $\cset{x_0,x_1,r_0,r'}$ is sufficient. Thus $\G\cong\C_2$, which contradicts the assumption.
    \end{enumerate}
    Hence we conclude that $f(u)=f(x_0)$ implies $u=x_0$.

    \noindent\textbf{Claim 2.} $f(a_0)\neq f(b_0)\neq f(b_1)\neq f(a_0)$.

    \noindent\textbf{Proof of Claim 2.} Suppose $f(a_0)=f(b_0)$. By Lemma~\ref{lem:t-morphism-source}, $f(x_0)\in f[\breve{R}[R[\breve{R}[b_0]]]]$. By Claim~1, the only possible case is that $f(r_0)=f(x_2)$, $f(r_1)=f(x_1)$ and $f(x_0)=f(u)$ for some $u\in X_L\cap\breve{R}[r_1]$. Then $f(r')\in f[R[x_2]]=\cset{f(a_0),f(x_1),f(x_2)}$. It is not hard to see that $f(r')=f(x_2)$ implies $\G\cong\C_1$, and $f(r')=f(x_1)$ contradicts to Claim~1. Thus $f(r')=f(a_0)=f(b_0)$, which entails $f(y_1)\in f[\breve{R}[r']]=\cset{f(a_0),f(x_2)}$. Since $\G\not\in\mathsf{TM}(\C_2)$, we have $f(y_1)=f(x_2)$ and $f(y_0)=f(x_1)$, which again contradicts to Claim~1. Thus $f(a_0)\neq f(b_0)$. 
    Note that $f(a_0)=f(b_1)$ implies $f(a_0)=f(b_0)$, we have $f(a_0)\neq f(b_1)$. 
    
    Suppose $f(b_0)=f(b_1)$. Since $f(y_1)\neq f(b_1)$, $f(y_1)=f(v)$ for some $v\in\breve{R}[b_1]\setminus\cset{b_1}$. Note that $v\in\breve{R}[a_0]$ and $f(a_0)\neq f(b_0)$, it is not hard to show that $f(y_0)=f(a_0)$, which entails that $f(x_0)\in f[\R^3[y_0]]$ and contradicts Claim~1.

    \noindent\textbf{Claim 3.} For all $n\in\omega$, $|Z_n|=2n+3$, where $Z_n=\cset{f(a_i):i\leq n}\cup\cset{f(b_i):i\leq n+1}$.

    \noindent\textbf{Proof of Claim 3.} The proof proceeds by induction on $n\in\omega$. The case $n=0$ follows from Claim~2 immediately. Let $n>0$. By induction hypothesis, it suffices to show $f(a_n),f(b_{n+1})\not\in f[Z_{n-1}]$ and $f(a_n)\neq f(b_{n+1})$. Since $\cset{f(a_{n-1}),f(b_{n-1})}\sub f[R[a_n]]\cap f[R[b_{n+1}]]$ and $\cset{f(a_{n-1}),f(b_{n-1})}\nsubseteq f[R[v]]$ for any $v\in Z_{n-1}$, we have $f(a_n),f(b_{n+1})\not\in f[Z_{n-1}]$. Note that $f(b_n)\in f[R[b_{n+1}]]\setminus f[R[a_n]]$, we see that $f(a_n)\neq f(b_{n+1})$.
    
    By Claim~3, $\G$ is infinite, which contradicts Lemma~\ref{lem:no-inf-Kripke-frame-S4tFi}. Hence $\Fr_r(\mathsf{Log}(\gf_I))=\mathsf{TM}(\C_2)$.
\end{proof}

\begin{lemma}\label{lem:same-frame-S4t}
    For all $I\sub\mathbb{Z}^+$, $\Fr(L)=\Fr(L_I)$.
\end{lemma}
\begin{proof}
    Since $L\not\in\NExt(\mathsf{S5}_t)$ and $\tup{\mathsf{Log}(\C_2),\mathsf{S5}_t}$ is a splitting pair in $\NExt(\ST)$, we have $\C_2\md L$ and so $\mathsf{TM}(\C_2)\sub\Fr_r(L)$. By Lemmas~\ref{lem:intersection-rooted-frame} and \ref{lem:no-nontrivial-image}, $\Fr_r(L_I)=\Fr_r(L)\cup\Fr_r(\mathsf{Log}(\gf_I))=\Fr_r(L)$. Hence $\Fr(L)=\Fr(L_I)$ for all $I\in\mathbb{Z}^+$.
\end{proof}

So far, we obtain a family of logics which share the same class of frames as $L$. It remains to show that the set $\cset{L_I:I\sub\mathbb{Z}^+}$ of tense logics is of the cardinality $2^{\aleph_0}$.

\begin{definition}
    Let $\phi_0:=\neg\mathsf{bd}_k[q_i/p_i]\wedge\bb\neg p$, where $p,q_0,\cdots,q_k\in\mathsf{Prop}$ are propositional variables which do not occur in $\phi_L$. Then we define
    \begin{itemize}
        \item $\phi_{x_0}:=\Delta^k\neg\phi_L\wedge\Delta^4_p\phi_0\wedge\nabla^3\neg\phi_0$, $\phi_{x_1}:=\bd\phi_{x_0}\wedge\neg\phi_{x_0}$ and $\phi_{x_2}:=\D\phi_{x_1}\wedge\neg\phi_{x_1}$;
        \item $\phi_{y_0}:=\Delta_p^7\phi_{x_0}\wedge\Delta^6\neg\phi_{x_0}$ and $\phi_{y_1}:=\D\phi_{y_0}\wedge\neg\phi_{y_0}$;
        \item $\phi_{a_0}:=\bd\phi_{x_2}\wedge\neg\phi_{x_2}$, $\phi_{b_0}:=\bd\phi_{y_1}\wedge\neg\phi_{y_1}$ and $\phi_{b_1}:=\D\phi_{b_0}\wedge\B\neg\phi_{a_0}\wedge\neg\phi_{b_0}$;
        \item $\phi_{AB}=\B(\phi_{b_0}\vee\phi_{b_1}\vee\D\bd\D\bd\phi_{x_0})$;
        \item for all $l\in\mathbb{Z}^+$, $\phi_{a_l}:=\phi_{AB}\wedge\D\phi_{a_{l-1}}\wedge\D\phi_{b_{l-1}}\wedge\B\neg\phi_{b_l}$ and $\phi_{b_{l+1}}:=\phi_{AB}\wedge\D\phi_{a_{l-1}}\wedge\D\phi_{b_{l}}\wedge\B\neg\phi_{a_l}$;
        \item for all $n\in\mathbb{Z}^+$, $\phi_{c_n}:=\neg\phi_{AB}\wedge\D\phi_{a_{n}}\wedge\B\neg\phi_{a_{n+1}}$.
    \end{itemize}
\end{definition}

\begin{lemma}\label{lem:sep-formulas-S4t}
    Let $U=A\cup B\cup\cset{x_0,x_1,x_2,y_0,y_1}$. For all $u\in U$ and $v\in X_I$, 
    \begin{enumerate}[(1)]
        \item $\gf_I,u\not\md\phi_u\to\nabla^k\phi_L$,
        \item for all valuation $V$ in $\gf_I$, $V(\phi_{x_0})\neq\ve$ implies $V(\phi_{AB})=A\cup B$ and $V(\phi_u)=\cset{u}$.
        \item $\gf_I,v\not\md\neg\phi_u$ implies $u=v$,
        \item $\gf_I\md\neg\phi_{c_j}$ for any $j\not\in I$.
    \end{enumerate}
\end{lemma}
\begin{proof}
    For (1), recall first that $\F_L,w_L\not\md\phi_L$. Then there exists a valuation $V':\mathsf{Prop}\to\mathcal{P}(X_L)$ such that $\F_L,V',w_L\md\neg\phi_L$. Note that $\gf_I$ is differentiated, there exists a valuation $V$ in $\gf_I$ such that the following conditions hold: (i) $V\rsto X_L=V'$, (ii) $V(p)=R[b_k]\cup\cset{x_0,x_1,x_2,y_0,y_1}$ and (iii) $V(q_i)=R[b_i]$ for all $i\leq k$. Let $\M=(\gf_I,V)$. Then $\M,b_{k+1}\md\phi_0$ and $\M,w_L\md\neg\phi_L$. Since $\gf_I$ is $k$-transitive, $\M,x_0\md\nabla^3\mathsf{bd}_k$ and there exists a $p$-path from $x_0$ to $b_{k+1}$ of length $4$, we see that $\M,x_0\md\phi_{x_0}$. Moreover, note that $\M,r_0\md\neg p$ and $V(\phi_0)\sub Y_I$, we have $\M,u\md\Delta^3\phi_0\vee\neg\nabla_p^4\phi_0$ for each point $u\in X_I\setminus\cset{x_0}$. Thus $x_0$ is the unique point satisfying $\phi_{x_0}$. By the construction the formulas $\phi_u$, the reader can now easily check that $V(\phi_{AB})=A\cup B$ and $\M,u\md\phi_u$ for all $u\in U$. Since $w_L\in\R^k[u]$, we have $\M,u\md\phi_u\wedge\Delta^k\neg\phi_L$.

    For (2), let $V$ be a valuation in $\gf_I$ and $\M=(\gf_I,V)$. Suppose $\M,v_0\md\phi_{x_0}$ for some $v_0\in X_I$. Then $\M,v_0\md\Delta_p^4\phi_0$ and $\M,v_1\md\phi_0$ for some $v_1\in X_I$. Since $\phi_0:=\neg\mathsf{bd}_k[q_i/p_i]\wedge\bb\neg p$, $R[v_1]$ contains a chain of length greater than $k$. Thus $v_1\in R[r_0]\setminus\cset{a_0,b_0,b_1,c_0}$ and so $\M,r_0\md\neg p$. Then $\M,w\md\Delta^3\phi_0\vee\neg\nabla_p^4\phi_0$ for each point $w\in X_I\setminus\cset{x_0}$, which entails that $V(\phi_{x_0})=\cset{x_0}$. Then clearly, $V(\phi_{AB})=A\cup B$ and $V(\phi_{w})=\cset{w}$ for all $w\in U$.

    For (3), take any $u\in U$ and $v\in X_I$. Suppose $\gf_I,v\not\md\neg\phi_u$. Then there exists a valuation $V$ in $\gf_I$ such that $\gf_I,V,v\md\phi_u$. Let $\M=(\gf_I,V)$. By the construction $\phi_u$, we always have $\md\phi_u\to\Delta^m\phi_{x_0}$ for some $m\in\omega$. Thus $V(\phi_{x_0})\neq\ve$. By (2), $V(\phi_{w})=\cset{w}$ for all $w\in U$, which entails $u=v$.

    For (4), take any $j\not\in I$. Suppose $\gf_I\not\md\neg\phi_{c_j}$. Then there exists $v\in X_I$ and a valuation $V$ in $\gf_I$ such that $\gf_I,V,v\md\phi_{c_j}$. Thus $V(\phi_{x_0})\neq\ve$. By (2),$V(\phi_{AB})=A\cup B$ and $V(\phi_w)=\cset{w}$ for $w\in\cset{a_{j},a_{j+1}}$. Since $j\not\in I$, we have $\breve{R}[a_j]\setminus(A\cup B\cup\breve{R}[a_{j+1}])=\ve$ and so $\M\md\neg\phi_{c_j}$, which is a contradiction. %
\end{proof}

\begin{lemma}\label{lem:pairwise-diff-S4t}
    For all $I,J\in\mathcal{P}(\mathbb{Z}^+)$, $I\neq J$ implies $L_I\neq L_J$.
\end{lemma}
\begin{proof}
    Take any distinct $I,J\in\mathcal{P}(\mathbb{Z}^+)$. Let $i\in I\setminus J$. It suffices to show that 
    \begin{center}
        $\phi_{c_i}\to\nabla^k\phi_L\in L_J\setminus L_I$.
    \end{center}
    By Lemma~\ref{lem:sep-formulas-S4t}(1), $\phi_{c_i}\to\nabla^k\phi_L\not\in\mathsf{Log}(\gf_I)\supseteq L_I$. By Lemma~\ref{lem:sep-formulas-S4t}(4), $\neg\phi_{c_i}\in\mathsf{Log}(\gf_J)$. Since $\phi_L\in L$, we see $\nabla^k\phi_L\in L$. Note that we may always assume that $\phi_L$ and $\phi_{c_i}$ contains no common variable, $\neg\phi_{c_i}\vee\nabla^k\phi_L\in L\cap\mathsf{Log}(\gf_J)$. Hence $\phi_{c_i}\to\nabla^k\phi_L\in L_J$.
\end{proof}

Note that $L\in\NExt(\ST)\setminus(\NExt(\mathsf{S5}_t)\cup\cset{\ST})$ is arbitrarily chosen, by Lemmas~\ref{lem:same-frame-S4t} and \ref{lem:pairwise-diff-S4t}, the following theorem holds:

\begin{theorem}\label{thm:degreeS4t-continuum}
    For all $L\in\NExt(\ST)\setminus(\NExt(\mathsf{S5}_t)\cup\cset{\ST})$, $\mathsf{deg}(L)=2^{\aleph_0}$.
\end{theorem}

Finally, we obtain our main results in this section: 
\begin{theorem}\label{thm:degreeS4t-main}
    Let $L\in\NExt(\ST)$. Then the following are equivalent:
    \begin{enumerate}[(1)]
        \item $L$ is an iterated splitting in $\NExt(\ST)$.
        \item $\mathsf{df}(L)=1$.
        \item $\mathsf{deg}(L)=1$.
        \item $\mathsf{df}(L)\neq 2^{\aleph_0}$.
        \item $\mathsf{deg}(L)\neq 2^{\aleph_0}$.
    \end{enumerate}    
\end{theorem}
\begin{proof}
    By Lemma~\ref{lem:deg=1}, (1) implies (2). Note that $|\NExt(\ST)|\leq 2^{\aleph_0}$, by Proposition~\ref{prop:deg-leq-df}, we see that (2) implies (3), as well as (4) implies (5). Clearly (3) implies (5), and (2) implies (4). It remains to note that (5) implies (1) follows from Theorem~\ref{thm:degreeS4t-continuum} and Lemma~\ref{lem:iterated-splitting-S4t}.
\end{proof}

It follows immediately that the following dichotomy theorem holds:
\begin{theorem}\label{thm:dichotomy-S4t}
    For all $L\in\NExt(\ST)$, $\mathsf{deg}(L)=\mathsf{df}(L)\in\cset{1,2^{\aleph_0}}$.
\end{theorem}

\begin{remark}\label{remark:better-iterated}
    As we can see, union-splittings in $\NExt(\ST)$ are still strictly Kripke-complete, while there exist strictly Kripke-complete logics in $\NExt(\ST)$ which are not union-splittings. Hence, so far, the notion of iterated splitting fits better with strictly Kripke-completeness, in the sense that for $L\in\cset{\mathsf{K},\mathsf{K}_t,\mathsf{K4}_t,\ST}$, the iterated splittings in $\NExt(L)$ are exactly the strictly Kripke-complete logics. 
\end{remark}

\section{Conclusions}

The present work contributes a series of results on the degree of Kripke-incompleteness in lattices of tense logics. We started with the lattice $\NExt(\mathsf{K}_t)$ of all tense logics. By giving a characterization of the degree of Kripke-incompleteness in $\NExt(\mathsf{K}_t)$, we proved the dichotomy theorem for $\NExt(\mathsf{K}_t)$, that is, $\mathsf{deg}_{\mathsf{K}_t}(L)\in\cset{1,2^{\aleph_0}}$ for all $L\in\NExt(\mathsf{K}_t)$. By the same method, we gave a characterization of the degree of Kripke-incompleteness in $\LT$ and proved the dichotomy theorem for $\LT$. Finally, we turned to Kripke-incompleteness in $\NExt(\ST)$. We showed that iterated splittings are strictly Kripke-complete and all other extensions of $\ST$ are of the degree $2^{\aleph_0}$. By showing that the degree of Kripke-incompleteness coincide with the degree of FMP in all the lattices mentioned above, we obtain also dichotomy theorem for the degree of FMP is these lattices.

We claim that we could obtain more results on the degree of Kripke-incompleteness in lattices of tense logics by the method given in this work. For example, consider the tense logic $\mathsf{K4D}^+_t=\mathsf{K4}_t\oplus\D\top$, which is the tense logic of serial frames. It is not hard to see that $\F_L\md\D\top$ implies $\gf_I\md\D\top$ for all $\gf_I$ defined in Section~\ref{subsec:degree-K4t}. Note that serial frames are closed under reflective unfolding, we claim that $\mathsf{deg}_{\mathsf{K4D}^+_t}(L)=2^{\aleph_0}$ for all proper extension $L$ of $\mathsf{K4D}^+_t$. Similarly, since $\F_L\md\mathsf{grz}^+$ implies $\gf_I\md\mathsf{grz}^+$ for all $\gf_I$ defined in Section~\ref{sec:DKI-S4t}, we claim that the dichotomy theorem holds for $\mathsf{Grz}^+=\mathsf{K}_t\oplus\mathsf{grz}^+$. Given that the frames $\F_I$ are frames for the bi-intuitionistic logic $\mathsf{biInt}$ and every finite bi-p-morphic image of $\F_I$ is a bi-p-morphic image of $\C_2$, we claim that the dichotomy theorem for the degree of FMP holds for $\mathsf{Ext}(\mathsf{biInt})$. These are left for future work.

There are still a lot of worth-studying future work and we outline a few topics here:

By Blok's characterization theorem, the union-splittings in $\NExt(\mathsf{K})$ are exactly the strictly Kripke-complete logics. Results obtained in this work showed that Blok's characterization theorem can be generalized to the lattices $\NExt(\mathsf{K}_t)$ and $\LT$. However, as we mentioned in Remark~\ref{remark:better-iterated}, every union-splitting in $\NExt(\ST)$ is strictly Kripke-complete while the inverse does not hold. Instead, iterated splittings fit perfectly with strictly Kripke-complete logics. So it is natural to ask: what is the relation between union-splittings, iterated splittings and strictly Kripke-complete logics in the lattices of tense logics? For example, is an iterated splitting always a union-splitting? Is it true that for all tense logic $L_0$ and $L$, $L$ is a union-splitting in $\NExt(L)$ implies $\mathsf{deg}_{L_0}(L)=1$? Is a strictly Kripke-complete tense logics always an iterated splitting?

As the reader might already notice, our method relies heavily on the reflective unfolding of Kripke frames. Recall that a tense logic $L$ is called finitely transitive if it is $n$-transitive for some $n\in\omega$. Our method used in this paper applies to only those lattices $\NExt(L_0)$ where $L_0$ is finitely transitive. An important future work is to study the degree of Kripke-incompleteness in lattices of finitely transitive tense logics, for example, $\NExt(\mathsf{S4.2}_t)$ and $\NExt(\mathsf{S4.3}_t)$. By \cite[Proposition~23]{Kracht1992}, there are infinite splittings in both of these lattices, which indicates that the case is quite different from the one for $\NExt(\ST)$. We believe these topics will require new methods and techniques.

As we have shown in Sections~\ref{sec:DKI-Kt-K4t} and \ref{sec:DKI-S4t}, the degree of FMP and the degree of Kripke-incompleteness coincide in the lattices $\K$, $\LT$ and $\NExt(\ST)$. The dichotomy theorem for the degree of FMP holds for all these three lattices. Bezhanishvili et al. \cite{Bezhanishvili.Bezhanishvili.ea2025} showed the anti-dichotomy theorem for the degree of FMP for $\NExt(\mathsf{K4})$ and $\NExt(\mathsf{S4})$. Here is a natural follow-up question: Is there any tense logic $L$ such that the anti-dichotomy theorem holds for $\NExt(L)$?

\bigskip

\noindent\textbf{Acknowledgement.} This work has been influenced by the \emph{Degree of Kripke-Incompleteness} project at the ILLC. The author is deeply grateful to Nick Bezhanishvili for his insightful and detailed comments, which significantly improved the manuscript. Thanks are also due to Tenyo Takahashi and Rodrigo Nicolau Almeida for the valuable discussions that helped shape the ideas presented here.

\bibliographystyle{acm}

\bibliography{../../0_Bibliography/References}

\end{document}